\documentclass[reqno, 12pt]{amsart}
\usepackage{amsmath,amsxtra,amssymb,latexsym, amscd,amsthm}
\usepackage{graphicx}
\usepackage[utf8]{inputenc}
\usepackage[mathscr]{euscript}
\usepackage{mathrsfs}
\usepackage[english]{babel}
\usepackage{color}
\usepackage{tikz}
\usepackage{hyperref}
\usepackage{enumerate}
\usepackage{graphicx} 
\usepackage[T1]{fontenc}

\newtheorem {theorem}{Theorem}[section]
\newtheorem {corollary}{Corollary}[section]
\newtheorem {proposition}{Proposition}[section]
\newtheorem {lemma}{Lemma}[section]
\newtheorem {example}{Example}[section]
\newtheorem {definition}{Definition}[section]
\newtheorem {remark}{Remark}[section]

\def\Limsup{\mathop{{\rm Lim}\,{\rm sup}}}

\def\ar{a\kern-.370em\raise.16ex\hbox{\char95\kern-0.53ex\char'47}\kern.05em}
\def\ees{{\accent"5E e}\kern-.385em\raise.2ex\hbox{\char'23}\kern-.08em}
\def\eex{{\accent"5E e}\kern-.470em\raise.3ex\hbox{\char'176}}
\def\AR{A\kern-.46em\raise.80ex\hbox{\char95\kern-0.53ex\char'47}\kern.13em}
\def\EES{{\accent"5E E}\kern-.5em\raise.8ex\hbox{\char'23 }}
\def\EEX{{\accent"5E E}\kern-.60em\raise.9ex\hbox{\char'176}\kern.1em}
\def\ow{o\kern-.42em\raise.82ex\hbox{
   \vrule width .12em height .0ex depth .075ex \kern-0.16em \char'56}\kern-.07em}
\def\OW{o\kern-.460em\raise1.36ex\hbox{
\vrule width .13em height .0ex depth .075ex \kern-0.16em
\char'56}\kern-.07em}
\def\UW{U\kern-.42em\raise1.36ex\hbox{
\vrule width .13em height .0ex depth .075ex \kern-0.16em \char'56}\kern-.07em}
\def\DD{D\kern-.7em\raise0.4ex\hbox{\char '55}\kern.33em}

\title[Pareto optimality under nondegeneracy at infinity]{Pareto optimality under nondegeneracy at infinity}

\author{Liguo Jiao}
\address[Liguo Jiao]{Academy for Advanced Interdisciplinary
Studies, Northeast Normal University, Changchun 130024,
Jilin Province, China}
\email{jiaolg356@nenu.edu.cn; hanchezi@163.com}

\author{Jae Hyoung Lee}
\address[Jae Hyoung Lee]{Department of Applied Mathematics, Pukyong National University, Busan 48513, Republic of Korea}
\email{mc7558@pknu.ac.kr; mc7558@naver.com}

\author[Ti\'{\^{e}}n-S\OW n Ph\d{a}m]{Ti\'{\^{e}}n-S\OW n Ph\d{a}m$^*$}
\address[Ti\ees n-S\ow n Ph\d{a}m]{Department of Mathematics, Dalat University, 1 Phu Dong Thien Vuong, Dalat, Vietnam}
\email{sonpt@dlu.edu.vn}

\thanks{}

\thanks{$^*$Corresponding author.}

\keywords{multiobjective optimization, (weakly/properly) efficient solutions, weak sharp minima, Newton polyhedra, nondegeneracy at infinity, semialgebraic}

\subjclass[2020]{Primary: 90C29~$\cdot$~90C23; secondary: 14P10}

\date{\today}

\begin{document}

\begin{abstract}
In this paper, we consider a broad class of multiobjective optimization problems, which are defined in a simple and elegant manner by using only two basic notions of Newton polyhedra and nondegeneracy at infinity associated with given data. For each problem in the class, its several basic properties are studied; in particular, necessary and sufficient conditions for the existence of (weakly/properly) efficient solutions as well as the H\"older-type global weak sharp minima are established.
\end{abstract}

\maketitle

\tableofcontents

\section{Introduction}\label{Section1}
This paper addresses some fundamental issues of global multiobjective optimization that revolve around the existence of (weakly/properly) efficient solutions and weak sharp minima. More precisely, let $f := (f_1, \ldots, f_m) \colon \mathbb{R}^n \to \mathbb{R}^m$ be a mapping and $\Omega$ be a nonempty closed subset of $\mathbb{R}^n.$ For any $z, z' \in \mathbb{R}^m,$ the notation $z \ge z'$ means that $z - z' \in \mathbb{R}^m_+,$ with $\mathbb{R}^m_+$ being the nonnegative orthant of $\mathbb{R}^m,$ i.e.,
\begin{eqnarray*}
\mathbb{R}_+^m &:=& \left\{z:=(z_1, \ldots, z_m) \in \mathbb{R}^m \mid z_i \geq 0, i = 1, \ldots, m \right\}.
\end{eqnarray*}
Then the multiobjective optimization problem
\begin{equation}\label{MP}
{\rm Min}_{\mathbb{R}_+^m} \left\{ f(x) \mid x \in \Omega\right\} \tag{MP}
\end{equation}
is the problem of finding all solutions. The notion in defining solutions for the problem~\eqref{MP} is that of {\em Pareto optimality,}
which was initially introduced by Pareto \cite{Pareto1989}. Recall that (see, for example, \cite{Ehrgott2005, Luc1989}) a point $\bar x \in \Omega$ is called an {\em efficient solution} (also known as {\em Pareto optimal solution}) to the problem~\eqref{MP} if it holds that
\begin{eqnarray*}
f(\Omega) \cap \left(f(\bar x) - \left(\mathbb{R}^m_+ \setminus \{0\} \right)\right) &=& \emptyset;
\end{eqnarray*}
a point $\bar x \in \Omega$ is called a {\em weakly efficient solution} (also known as {\em weak Pareto optimal solution}) to the problem~\eqref{MP} if it holds that
\begin{eqnarray*}
f(\Omega) \cap \left(f(\bar x) - {\rm int}\mathbb{R}^m_+ \right) &=& \emptyset.
\end{eqnarray*}
Recall also from \cite{Geoffrion1968} that a point $\bar x \in \Omega$ is said to be a {\em properly efficient solution (in Geoffrion's sense)} to the problem~\eqref{MP} if it is an efficient solution and there is a real number $M > 0$ such that for each {\mbox{$i\in\{1,\ldots,m\}$}} and each $x\in\Omega$ satisfying $f_i(x) < f_i(\bar x)$ there exists at least one $j\in\{1,\ldots,m\}$ with $f_j(\bar x) < f_j(x)$ and
\begin{eqnarray*}
\frac{f_i(\bar x) - f_i(x)}{f_j(x) - f_j(\bar x)}&\le& M.
\end{eqnarray*}
When $\bar x$ is a (weakly/properly) efficient solution, $f(\bar x)$ is called a {\em $($weakly/properly$)$ efficient value.} 

We are concerned with the question whether the problem~\eqref{MP} has a (weakly/properly) efficient solution.
The answer to this question can be trivial if, for example, the objective mapping $f$ is continuous and the constraint set $\Omega$ is compact, or if the restriction $f|_\Omega$ has no empty sections (definitions will be given later). In the former case, the answer is positive by the well known Weierstrass theorem. In the latter case, no solution can exist. The situation, however, becomes more subtle when $\Omega$ is unbounded.
To address this challenging unbounded case, significant progress has been made by tailoring existence conditions to specific classes of problems; see, e.g., \cite{Bao2010, Borwein1983, Ha2006, Ha2010, Kim2021, Kim2019}. 
For instance, when $f$ is a linear mapping and $\Omega$ is a polyhedral set, some existence theorems have been established under asymptotic cone conditions; see, e.g., \cite[Theorem 4.2.1, Corollary 4.2.5 and Exercise 4.4.21]{Luc2016}. 
In contrast, when $f$ is a locally Lipschitz mapping and the restriction $f|_\Omega$ has a nonempty bounded section, several existence theorems have been elegantly derived under Palais--Smale conditions; see \cite{Bao2010, Ha2006, Kim2021, Kim2019}.

We are also interested in the {\em existence of weak sharp minima}. In scalar optimization, the notion of a {\it sharp minimum} was introduced by Polyak~\cite{Polyak1987}. Later, Ferris \cite{Ferris1988} (see also \cite{BurkeFerris1993}) proposed the terminology {\it weak sharp minima} to include the possibility of non-unique solution set. A weak sharp minimum is a property of the solution set of an optimization problem that  says that the objective function grows at least linearly with the distance from any feasible point to that set. It has been well recognized that weak sharp minima play an important role in the stability/sensitivity analysis of optimization problems and the convergence analysis of some numerical algorithms.
In the last decade, there has been increasing interest in weak sharp minima property in the context of multiobjective optimization; see, e.g., \cite{Deng2004, Hu2024, Liu2009, PHAMTS2018-1, Studniarski2007, YangYen2010, ZhengNg2014, ZhengYang2008} and the references therein. In particular, it is shown in \cite{Deng2004, Hu2024, YangYen2010, ZhengYang2008} that the set of (weakly) efficient solutions of a linear multiobjective optimization problem has the global weak sharp minimum property. In contrast, the set of weakly efficient solutions of a multiobjective optimization problem with polynomial data, which satisfies a regularity condition at infinity, has the H\"older type global weak sharp minimum property; see~\cite{PHAMTS2018-1}.

In this paper, we focus on a broad class of multiobjective optimization problems with semialgebraic data
in which the existence of (weakly/properly) efficient solutions is ensured and weak sharp minima hold. Our approach is mainly based on tools of variational analysis and semialgebraic geometry.

\subsubsection*{Organization and main contributions}
In Section~\ref{Section2}, some definitions and preliminary results from variational analysis and semialgebraic geometry are recalled.

In Section~\ref{Section3}, the concepts of Newton polyhedra and nondegeneracy at infinity are given. The latter one extends the concepts defined by 
Kouchnirenko~\cite{Kouchnirenko1976} and Khovanskii~\cite{Khovanskii1978} for complex polynomial mappings which play an important role in Singularity Theory (see, for example, \cite{Arnold2012, Oka1997}). We show that, in the space of polynomial mappings with fixed Newton polyhedra, the class of polynomial mappings, whose restrictions on a closed semialgebraic set are nondegenerate at infinity, forms an open set.

In Section~\ref{Section4}, we prove a key result in later developments which says that if $f$ is convenient and $f|_\Omega$ is nondegenerate at infinity, then for every continuous selection of the component polynomials of $f,$ its restriction to $\Omega$ has good behavior at infinity. This leads to a continuity result for the solution set defined by a system of polynomial inequalities.

In Section~\ref{Section5}, we consider the problem~\eqref{MP} with $f$ being convenient and $f|_\Omega$ being nondegenerate at infinity. We first present necessary and sufficient conditions for the existence of (weakly/properly) efficient solutions. Then we establish the H\"older-type global weak sharp minimum property as well as the upper semicontinuity of sublevel sets.

Finally, in Section~\ref{Section6}, some examples and remarks are given.


\section{Preliminaries} \label{Section2}

We begin by giving some necessary definitions and notational conventions. Let $\mathbb{R}^n$ denote the Euclidean space of dimension $n.$ The corresponding inner product (resp., norm) in $\mathbb{R}^n$ is defined by $\langle x, y \rangle$ for any $x, y \in \mathbb{R}^n$ (resp., $\| x \| := \sqrt{\langle x, x \rangle}$ for any $x \in \mathbb{R}^n$). 
The closed ball (resp., the sphere) centered at $\bar{x} \in \mathbb{R}^n$ of radius $r$ will be denoted by $\mathbb{B}_{r}(\bar{x})$ (resp., $\mathbb{S}_{r}(\bar{x})$); when $\bar{x}$ is the origin of $\mathbb{R}^n$ we write $\mathbb{B}_{r}$ (resp., $\mathbb{S}_r$) instead of $\mathbb{B}_{r}(\bar{x})$ (resp., $\mathbb{S}_{r}(\bar{x})).$ We will adopt the convention that $\inf \emptyset = \infty$ and $\sup \emptyset = -\infty.$
The closure, interior, and boundary of a set ${S}$ are denoted by $\mathrm{cl} {{S}},$ $\mathrm{int} {S},$ and $\mathrm{bd} {S},$ respectively.
A set $S \subset \mathbb{R}^n$ is said to be {\em locally closed} if for each $x \in S$ there exists a neighborhood $U \subset \mathbb{R}^n$ of ${x}$ such that $S \cap U$ is a closed set.

For an extended real-valued function $f \colon \mathbb{R}^n \rightarrow \mathbb{R} \cup \{\infty\},$ we denote its {\em effective domain} and {\em epigraph} by, respectively,
\begin{eqnarray*}
\mathrm{dom} f &:=& \{ x \in \mathbb{R}^n \ | \ f(x) < \infty\},\\
\mathrm{epi} f &:=& \{ (x, r) \in \mathbb{R}^n \times \mathbb{R} \ | \ f(x) \le r \}.
\end{eqnarray*}
The function $f$ is {\em lower semicontinuous} if for every $x \in \mathbb{R}^n$ the inequality $\liminf_{x' \to {x}} f(x') \ge f({x})$ holds, or equivalently, the epigraph set $\mathrm{epi} f$ is closed in $\mathbb{R}^n \times \mathbb{R} .$

The {\em indicator function} of a set ${S} \subset \mathbb{R}^n,$ denoted by $\delta_{S},$ is defined by
\begin{eqnarray*}
\delta_{S}(x) &:=&
\begin{cases}
0 & \textrm{ if } x \in {S}, \\
\infty & \textrm{ otherwise.}
\end{cases}
\end{eqnarray*}
By definition, $\delta_{S}$ is lower semicontinuous if and only if ${S}$ is closed.

The {\em Painlev\'e--Kuratowski outer limit} of a set-valued mapping $\Phi \colon \mathbb{R}^n \rightrightarrows \mathbb{R}^m$ at $x$ is defined by
\begin{eqnarray*}
\Limsup_{x' \to {x}} \Phi(x') &:=& \left\{y \in \mathbb{R}^m \mid \exists x_k \to {x},\ \exists y_k \in \Phi(x_k),\ y_k \to y\right\}.
\end{eqnarray*}

\subsection{Normals and subdifferentials}

Here we recall some concepts and properties of normal cones to sets and subdifferentials of real-valued functions, which can be found in~\cite{Mordukhovich2018, Rockafellar1998}.

\begin{definition}{\rm
Consider a set ${S}\subset\mathbb{R}^n$ and a point ${x} \in {S}.$
\begin{enumerate}[\upshape (i)]
\item The {\em regular normal cone} (also known as the {\em prenormal} or {\em Fr\'echet normal cone}) $\widehat{N}({x}; {S})$ to
${S}$ at ${x}$ consists of all vectors $v\in\mathbb{R}^n$ satisfying
\begin{eqnarray*}
\limsup_{x' \xrightarrow{S} {x}} \frac{\langle v , x' - {x} \rangle}{\|x' -{x}\|} & \leq & 0,
\end{eqnarray*}
where $x' \xrightarrow{S} {x}$ means that $x' \rightarrow {x} $ with $x' \in S.$

\item The {\em limiting normal cone} (also known as the {\em basic} or {\em Mordukhovich normal cone}) $N({x}; {S})$ to ${S}$ at ${x}$ consists of all vectors $v \in \mathbb{R}^n$ such that there are sequences $x_k \to {x}$ with $x_k \in {S}$ and $v_k \rightarrow v$ with $v_k \in \widehat N(x_k; {S}),$ or in other words,
\begin{eqnarray*}
{N}({x}; S) &:=& \Limsup_{x' \xrightarrow{S} {x}}\widehat{N}({x}'; S).
\end{eqnarray*}
\end{enumerate}

If $x \not \in {S},$ we put $\widehat{N}(x; {S}) := \emptyset$ and ${N}(x; {S}) := \emptyset.$
}\end{definition}

\begin{remark}{\rm
If ${S}$ is a manifold of class $C^1,$ then for every point $x \in {S},$ the normal cones $\widehat{N}({x}; {S})$ and $N({x}; {S})$ are equal to the normal space to ${S}$ at ${x}$ in the sense of differential geometry; see \cite[Example~6.8]{Rockafellar1998}. In particular, for all $R > 0$ and all $x \in \mathbb{S}_R,$ we have $N(x; \mathbb{S}_R) = \{\lambda x \mid \lambda \in \mathbb{R}\}.$
}\end{remark}

Functional counterparts of normal cones are subdifferentials.
\begin{definition}{\rm
Consider a function $f\colon\mathbb{R}^n \to \mathbb{R} \cup \{\infty\}$ and a point ${x} \in \mathrm{dom} f.$ The {\em limiting} and {\em singular subdifferentials} of $f$ at ${x}$ are defined by, respectively,
\begin{eqnarray*}
\partial f({x}) &:=& \big\{v\in\mathbb{R}^n\;\big|\;(v,-1)\in N\big(({x}, f({x}));\mathrm{epi} f\big)\big\}, \\
\partial^\infty f({x}) &:=& \big\{v\in\mathbb{R}^n\;\big|\;(v,0)\in N\big(({x}, f({x}));\mathrm{epi} f \big)\big\}.
\end{eqnarray*}
}\end{definition}

In \cite{Mordukhovich2018, Rockafellar1998} the reader can find equivalent analytic descriptions of the limiting subdifferential $\partial f({x})$ and comprehensive studies of it and related constructions. For convex $f,$ this subdifferential coincides with the convex subdifferential. Furthermore, if the function $f$ is of class $C^1,$ then $\partial f({x}) = \{\nabla f({x})\}.$ The singular subdifferential $\partial^\infty f({x})$ plays an entirely different role, it detects horizontal ``normal'' to the epigraph, and it plays a decisive role in qualification conditions of subdifferential calculus.

The following facts are well known (see, for example, \cite{Mordukhovich2018, Rockafellar1998}).

\begin{lemma}\label{Lemma23}
For any locally closed set $S \subset \mathbb{R}^n$ and point ${x} \in S,$ we have
\begin{eqnarray*}
\partial \delta_S({x}) &=& \partial^\infty \delta_S({x}) \ = \ N ({x}; S).
\end{eqnarray*}
\end{lemma}

\begin{lemma}\label{FermatRule}
Let $S \subset \mathbb{R}^n$ be a locally closed set around $\overline{x} \in S$ and let $f \colon \mathbb{R}^n \rightarrow {\mathbb{R}}$ be a Lipschitz function in some neighborhood of $\overline{x}.$ If $\overline{x}$ is a local minimizer of $f$ on $S$ then
\begin{eqnarray*}
0 & \in & \partial f(\overline{x}) + N(\overline{x}; S).
\end{eqnarray*}
\end{lemma}

Recall that a function $h \colon \mathbb{R}^n \to \mathbb{R}$ is said to be a {\em continuous selection} of a given collection of continuous functions $f_1, \ldots, f_m$ on $\mathbb{R}^n$ if $h$ is continuous and $h(x) \in \{f_1(x), \ldots, f_m(x)\}$ for every $x \in \mathbb{R}^n.$

\begin{lemma}[{see \cite[Theorem~4.3.1]{Scholtes2012}}]\label{Lemma24}
Let $h \colon \mathbb{R}^n \to \mathbb{R}$ be a continuous selection of a collection of $C^1$-functions $f_i \colon \mathbb{R}^n \to \mathbb{R}, i = 1, \ldots, m.$ Then for any $x \in \mathbb{R}^n,$ we have
\begin{eqnarray*}
\partial h ({x}) &\subset& \left\{ \sum_{i\in I({x})} \lambda_i \nabla f_i({x})\;\Big|\; \lambda_i\ge 0 \ \  \textrm{ and } \ \ \sum_{i\in I({x})}\lambda_i=1\right\},
\end{eqnarray*}
where $I({x}) := \{i \mid f_i({x}) = h({x})\}.$
\end{lemma}

\subsection{Semialgebraic geometry}

Now, we recall some notions and results of semialgebraic geometry, which can be found in~\cite{Bochnak1998} and \cite[Chapter~1]{HaHV2017}.

\begin{definition}{\rm
\begin{enumerate}[{\rm (i)}]
\item A subset $S$ of $\mathbb{R}^n$ is called {\em semialgebraic} if it is a finite union of sets of the form
$$\{x \in \mathbb{R}^n \mid f_i(x) = 0, \ i = 1, \ldots, p; f_i(x) > 0, \ i = p + 1, \ldots, q\},$$
where all $f_{i}$ are polynomials.
In other words, $S$ is a union of finitely many sets, each defined by finitely many polynomial equalities and inequalities.

\item A mapping $f \colon S \rightarrow ({\mathbb{R} \cup \{\infty\})^m}$ is said to be {\em semialgebraic} if its graph
\begin{equation*}
\{ (x, y) \in S \times \mathbb{R}^m \mid y =f(x) \}
\end{equation*}
is a semialgebraic subset of $\mathbb{R}^n \times \mathbb{R}^m.$
\end{enumerate}
}\end{definition}

A major fact concerning the class of semialgebraic sets is its stability under linear projections.

\begin{lemma}[\textnormal{Tarski--Seidenberg theorem}] \label{TarskiSeidenbergTheorem}
The image of any semialgebraic set $S \subset \mathbb{R}^n$ under a projection to any linear subspace of $\mathbb{R}^n$ is a semialgebraic set.
\end{lemma}

\begin{remark}\label{TarskiSeidenbergRemark}{\rm
As an immediate consequence of the Tarski--Seidenberg Theorem, we get semialgebraicity of any set $\{ x \in A \mid \exists y \in B,(x, y) \in C \},$ provided that $A ,B,$ and $C$ are semialgebraic sets in the corresponding spaces. Also, $\{ x \in A \mid  \forall y \in B,(x, y) \in C \}$ is a semialgebraic set as its complement is the union of the complement of $A$ and the set $\{ x \in A \mid  \exists y \in B,(x, y) \not\in C \}.$ Thus, if we have a finite collection of semialgebraic sets, then any set obtained from them with the help of a finite chain of quantifiers is also semialgebraic.
}\end{remark}

The following well known lemmas will be of great importance for us (see, for example, \cite{Bochnak1998, HaHV2017}).

\begin{lemma} [\textnormal{monotonicity lemma}] \label{MonotonicityLemma}
Let $f \colon (a, b) \rightarrow \mathbb{R}$ be a semialgebraic function. Then there are finitely many points $a =: t_0 < t_1 < \cdots < t_p := b$ such that for each $i = 0, \ldots, p - 1,$ the restriction of $f$ to the interval $(t_i, t_{i + 1})$ is analytic$,$ and either constant$,$ or strictly increasing or strictly decreasing.
\end{lemma}

\begin{lemma}[\textnormal{growth dichotomy lemma}] \label{GrowthDichotomyLemma}
Let $f \colon (0, \epsilon) \rightarrow {\mathbb R}$ be a semialgebraic function with $f(t) \ne 0$ for all $t \in (0, \epsilon).$ Then there exist constants $a \ne 0$ and $r \in {\mathbb Q}$ such that $f(t) = at^{r} + o(t^{r})$ as $t \to 0^+.$
\end{lemma}

\begin{lemma}[curve selection lemma] \label{CurveSelectionLemma1}
Let $S \subset \mathbb{R}^n$ be a semialgebraic set and $x^* \in \mathbb{R}^n$ be a nonisolated point of $\overline{S}.$
Then there exists a nonconstant analytic semialgebraic curve $\varphi\colon (-\epsilon, \epsilon) \to {\mathbb R}^n$ with $\varphi(0) = x^*$ and with $\varphi(t) \in S$ for $t \in (0, \epsilon).$
\end{lemma}

\begin{lemma}[\textnormal{curve selection lemma at infinity}]\label{CurveSelectionLemmaAtInfinity}
Let $S\subset \mathbb{R}^n$ be a semialgebraic set$,$ and let $f := (f_1, \ldots,f_m) \colon\mathbb{R}^n \to \mathbb{R}^m$ be a semialgebraic mapping. Assume that there exists a sequence $\{x_k\}_{{k} \ge 1} \subset S$ such that $\lim_{k \to \infty} \| x_k\| = \infty$ and $\lim_{k \to \infty} f(x_k)= y \in( \mathbb{R} \cup \{\pm \infty\})^m.$ Then there exists an analytic semialgebraic curve $\phi \colon (0, \epsilon) \to \mathbb{R}^n$ lying in $S$ such that $\lim_{t\to 0^+} \|\phi(t)\| = \infty$ and $\lim_{t \to 0^+} f(\phi(t)) = y.$
\end{lemma}

We also need the following facts.
\begin{lemma}[{see \cite[Lemma~2.10]{Drusvyatskiy2015}}] \label{ChainRule}
Consider a lower semicontinuous, semialgebraic function $f \colon \mathbb{R}^n \rightarrow \mathbb{R} \cup \{\infty\}$ and a semialgebraic curve $\phi \colon (0, \epsilon) \rightarrow {\rm dom}f.$
Then for all $t > 0$ sufficiently small$,$ the mappings $\phi$ and $f \circ \phi$ are analytic at ${{t}}$ and satisfy
\begin{eqnarray*}
\left \langle v(t), \frac{d \phi(t)}{dt} \right \rangle &=& \frac{d}{dt}(f\circ\phi)(t) \quad \textrm{ for all } \quad v(t) \in \partial f (\phi({{t}})).
\end{eqnarray*}
\end{lemma}

\begin{corollary}\label{Corollary214}
Let $S\subset \mathbb{R}^n$ be an unbounded closed semialgebraic set. Then there exists a real number $R > 0$ such that for all $t \ge R$
the set $S \cap \mathbb{S}_t$ is nonempty, and we have for all $x \in S \cap \mathbb{S}_t$ that
\begin{eqnarray*}
N(x; S)\cap \big(-N(x; \mathbb{S}_t)\big) &=& \{0\},
\end{eqnarray*}
in particular,
\begin{eqnarray*}
N(x; S \cap \mathbb{S}_t) &\subset& N(x; S) + N(x; \mathbb{S}_t).
\end{eqnarray*}
\end{corollary}
\begin{proof}
By \cite[Theorem~2.16]{Mordukhovich2018}, it suffices to show the first statement.
Assume for contradiction that there is an analytic semialgebraic curve $(\phi, \mu) \colon (R, \infty) \rightarrow \mathbb{R}^n \times \mathbb{R}$ such that for all $t > R,$
\begin{eqnarray*}
\mu(t) \ne 0, \quad \|\phi({{t}})\| = {{t}}, \quad \textrm{ and } \quad -\mu({{t}}) \phi({{t}}) \in N(\phi({{t}}; S).
\end{eqnarray*}
Applying Lemma~\ref{ChainRule} and increasing $R$ (if necessary), we have for any ${{t}} > R,$
\begin{eqnarray*}
0 &=& \langle -\mu({{t}}) \phi({{t}}),
\frac{d \phi(t)}{dt} \rangle \ = \
- \mu({{t}}) \frac{1}{2} \frac{d}{d {{t}}} \|\phi({{t}})\|^2 \ = \ - \mu({{t}}){{t}},
\end{eqnarray*}
which yields $\mu(t) = 0,$ a contradiction.
\end{proof}

\section{Newton polyhedra and the nondegeneracy at infinity} \label{Section3}

In this section, the concepts of Newton polyhedra and nondegeneracy at infinity are given. Let us start with some notation and definitions.

\subsection{Newton polyhedra}
We denote by $\mathbb{Z}_+$ and $\mathbb{R}_+$, respectively, the set of nonnegative integers and the set of nonnegative real numbers.
For a nonempty subset $J$ of $\{1, \ldots, n\},$ we define
\begin{eqnarray*}
\mathbb{R}^J &:=& \{x \in \mathbb{R}^n \mid x_j = 0 \ \textrm{ for all } \ j \not \in J\}, \\
\mathbb{(R^*)}^J &:=& \{x \in \mathbb{R}^n \mid x_j = 0 \ \textrm{ if and only if } \ j \not \in J\},
\end{eqnarray*}
$\mathbb{R}^J_+ := \mathbb{R}^J \cap \mathbb{R}^n_+$ and $\mathbb{(R_+^*)}^J  := \mathbb{(R^*)}^J  \cap \mathbb{R}^n_+.$ If $\kappa = (\kappa_1, \ldots, \kappa_n) \in \mathbb{Z}_+^n,$ we denote by $x^\kappa$ the monomial $x_1^{\kappa_1} \cdots x_n^{\kappa_n}.$

A subset $\Gamma \subset {\mathbb R}^n_+$ is a {\em Newton polyhedron} if there exists a finite subset $A \subset {\mathbb Z}^n_+$  such that $\Gamma$ is the convex hull in ${\mathbb R}^n$ of $A.$ We say that $\Gamma$ is the {\em Newton polyhedron determined by $A$} and write $\Gamma = \Gamma(A).$

Given a Newton polyhedron $\Gamma$ and a vector $q \in {\mathbb R}^n,$ we let
\begin{eqnarray*}
d(q, \Gamma) &:=& \min \{\langle q, \kappa \rangle \mid \kappa \in \Gamma\}, \\
\Delta(q, \Gamma) &:=& \{\kappa \in \Gamma \mid \langle q, \kappa \rangle = d(q, \Gamma) \}.
\end{eqnarray*}
By definition, for each nonzero vector $q \in \mathbb{R}^n,$
$\Delta(q, \Gamma) $ is a closed face of $\Gamma.$ Conversely,
if $\Delta$ is a closed face of $\Gamma$, then there exists a
nonzero vector $q\in \mathbb{R}^n$ such that $\Delta = \Delta(q, \Gamma),$ where we
can in fact assume that $q\in\mathbb{Q}^n$ since $\Gamma$ is an integer polyhedron.
The {\em dimension} of a face $\Delta$ is the minimum of the
dimensions of the affine subspaces containing $\Delta.$
The faces of $\Gamma$ of dimension $0$ are the {\em vertices} of $\Gamma.$

Let $f \colon \mathbb{R}^n \to \mathbb{R}$ be a polynomial function. Suppose that $f$ is written as $f = \sum_{\kappa} c_\kappa x^\kappa.$ The {\em support} of $f,$ denoted by $\mathrm{supp}(f),$ is the set of $\kappa  \in \mathbb{Z}_+^n$ such that $c_\kappa \ne 0.$
The {\em Newton polyhedron (at infinity)} of $f$, denoted by $\Gamma(f),$ is the convex hull in $\mathbb{R}^n$ of the
set\footnote{We do not include the origin $0 \in \mathbb{R}^n$ in the convex hull of  $\mathrm{supp}(f),$ as is usual. It turns out that we do not need the origin to obtain our results, as we shall see. Since taking the origin into account is really unnecessary in our context, we omit it.}
 $\mathrm{supp}(f),$ i.e.,  $\Gamma(f)=\Gamma(\mathrm{supp}(f)).$
For each (closed) face $\Delta$ of $\Gamma(f),$  we will denote
\begin{eqnarray*}
f_\Delta(x) &:=& \sum_{\kappa \in \Delta \cap \mathbb{Z}_+^n} c_\kappa x^\kappa.
\end{eqnarray*}

\begin{remark}{\rm
Let $\Delta := \Delta(q, \Gamma(f))$ for some nonzero vector $q := (q_1, \ldots, q_n) \in \mathbb{R}^n.$ By definition, $f_\Delta$ is a weighted homogeneous polynomial of type $(q, d := d(q, \Gamma(f))),$ i.e., we have for all $t > 0$ and all $x \in \mathbb{R}^n,$
\begin{eqnarray*}
f_\Delta(t^{q_1} x_1,  \ldots, t^{q_n} x_n) &=& t^d f_\Delta(x_1, \ldots, x_n).
\end{eqnarray*}
}\end{remark}

\subsection{The cone $C(\infty; \Omega)$}
For any $x := (x_1, \ldots, x_n)$ and $v := (v_1, \ldots, v_n)$ we write $x \odot v := (x_1 v_1, \ldots, x_n v_n),$ which is a vector in $\mathbb{R}^n.$
For an unbounded subset $\Omega$ of $\mathbb{R}^n$, we let (see \cite{Jiao2024, PHAMTS2025-1})
\begin{eqnarray*}
C(\infty; \Omega) &:=& \Limsup_{x \xrightarrow{\Omega} \infty} x \odot \widehat{N}(x; {\Omega}).
\end{eqnarray*}
Each $\xi \in C(\infty; \Omega)$ is presented as follows: there are sequences $x_k \in \Omega$ and $v_k \in \widehat{N}(x_k; \Omega)$ such that $\|x_k\| \to \infty$ and $x_k \odot v_k \to \xi$ as $k \to \infty.$ It is easy to see that $C(\infty; \Omega)$ is a closed cone in $\mathbb{R}^n$ (not necessarily convex). Moreover, we have:

\begin{lemma}[{see \cite[Proposition~1]{PHAMTS2025-1}}] \label{MD32}
The following equality holds$:$
\begin{eqnarray*}
C(\infty; \Omega) &=& \Limsup_{x \xrightarrow{\Omega} \infty} x \odot N(x; {\Omega}).
\end{eqnarray*}
\end{lemma}

\subsection{The nondegeneracy at infinity}

Let $f := (f_1, \ldots, f_m) \colon \mathbb{R}^n \to \mathbb{R}^m$ be a polynomial mapping and $\Omega \subset \mathbb{R}^n$ be an unbounded, closed and semialgebraic set. The following definition is inspired by those of Kouchnirenko \cite{Kouchnirenko1976} and Khovanskii \cite{Khovanskii1978} (see also \cite{HaHV2017}).

\begin{definition}\label{nondegeneracy}{\rm
The restriction $f|_\Omega$ is called {\em nondegenerate at infinity} if for any vector\footnote{The number of vectors $q \in \mathbb{R}^n$ is infinite; however, there exists a finite number of faces $\Delta_i$ of the Newton polyhedron $\Gamma(f_i).$}  $q \in \mathbb{R}^n$ with $\min_{j = 1, \ldots, n} q_j < 0$ there are no nonempty sets $I \subset \{1, \ldots, m\}$ and $J \subset \{1, \ldots, n\},$ a vector $c \in (\mathbb{R}^*_+)^I$ and
a point $x^* \in \mathbb{(R^*)}^J$ such that the following conditions hold:
\begin{enumerate}[\upshape (c1)]
\item $d(q, \Gamma(f_i)) < 0$ and  $\Delta_i := \Delta(q, \Gamma(f_i)) \subset \mathbb{R}^J$ for all $i \in I;$
\item $f_{i, \Delta_i}(x^*) = 0$ for all $i \in I;$
\item $0 \in \sum_{i \in I} c_i x^* \odot \nabla f_{i, \Delta_i} (x^*) + C(\infty; \Omega) \cap \mathbb{R}^J.$
\end{enumerate}
}\end{definition}

When $\Omega$ is the whole space $\mathbb{R}^n,$ we say ``the mapping $f$ is nondegenerate at infinity'' instead of
``the restriction $f|_\Omega$ is nondegenerate at infinity''. In this case, $C(\infty; \Omega) = \{0\}$ and the condition (c3) is equivalent to the fact that
\begin{eqnarray*}
\sum_{i \in I}c_i \nabla f_{i, \Delta_i}(x^*) &=& 0,
\end{eqnarray*}
and hence the nondegeneracy at infinity in the above definition is weaker than the one in \cite[Definition 3.8]{HaHV2017}, because we require here that $c_i > 0$ for all $i \in I.$

\subsection{The openness of the nondegeneracy condition at infinity}

We close this section by proving that in the space of polynomial mappings with fixed Newton polyhedra, the set of polynomial mappings, whose restrictions on a closed and semialgebraic set are nondegenerate at infinity, is an open set. To this end, we start with the following simple lemma.

\begin{lemma} \label{PolyhedralLemma}
Let $\Gamma \subset \mathbb{R}^n$ be a Newton polyhedron and $q \in \mathbb{R}^n.$ Then for all $\widetilde{q} \in \mathbb{R}^n$ close to $q$ we have
\begin{eqnarray*}
\Delta(\widetilde q, \Gamma) &= & \Delta(\widetilde q, \Delta(q, \Gamma)).
\end{eqnarray*}
\end{lemma}
\begin{proof}
We first show that for all $\widetilde{q} \in \mathbb{R}^n$ close to $q,$
\begin{eqnarray} \label{Eqn1}
\Delta(\widetilde q, \Gamma) & \subset&\Delta(q, \Gamma).
\end{eqnarray}
Indeed, if this were not true, there would exist a sequence $q_k$ tending to $q$ such that $\Delta(q_k, \Gamma) \setminus \Delta(q, \Gamma) \ne \emptyset.$ It is well known that $\Delta(q_k, \Gamma)$ and $\Delta(q, \Gamma)$ are convex hulls of some subsets of vertices of $\Gamma$ and that
the number of vertices of $\Gamma$ is finite. Passing to a subsequence if necessary, we may assume that there exists a vertex $\kappa \in \Delta(q_k, \Gamma) \setminus \Delta(q, \Gamma)$ for all $k.$ Then the following relations hold
$$\langle q_k, \kappa \rangle= d(q_k, \Gamma), \quad \langle q, \kappa \rangle>d(q, \Gamma)
\quad \textrm{ and } \quad \langle q_k, \kappa\rangle \to \langle q, \kappa\rangle \ \textrm{ as } \ k \to \infty.$$
On the other hand, by definition, it is not hard to check that
\begin{eqnarray*}
|d(q_k, \Gamma) - d(q, \Gamma) | &\leq & \|q_k - q\| \sup_{\kappa' \in \Gamma} \|\kappa'\|.
\end{eqnarray*}
Consequently, $d(q_k, \Gamma) \to d(q, \Gamma)$ as $k \to \infty,$ which contradicts the relations mentioned above.

Next, let $\widetilde{q}$ satisfy the inclusion~\eqref{Eqn1} and we will prove the desired equality. For simplicity, write
$\Delta := \Delta(q, \Gamma).$ Take any $\kappa \in \Delta(\widetilde q, \Gamma).$ We have $\kappa \in \Delta$ (by the inclusion~\eqref{Eqn1}) and
$\langle \widetilde q, \kappa \rangle= d(\widetilde q, \Gamma) \le d(\widetilde q, \Delta)$ (since $\Gamma \supset \Delta).$
By definition, $\kappa \in \Delta(\widetilde q, \Delta).$ Therefore,
\begin{eqnarray*}
\Delta(\widetilde q, \Gamma) & \subset & \Delta(\widetilde q, \Delta).
\end{eqnarray*}

For the converse inclusion, take any $\kappa \in \Delta(\widetilde q, \Delta).$ We have
\begin{eqnarray*}
\langle \widetilde q, \kappa \rangle&=& d(\widetilde q, \Delta) \ \le \ d(\widetilde q, \Delta(\widetilde q, \Gamma)) \ = \ d(\widetilde q, \Gamma),
\end{eqnarray*}
where the inequality follows from the inclusion~\eqref{Eqn1}. Consequently, $\kappa \in \Delta(\widetilde q, \Gamma).$
Therefore, $\Delta(\widetilde q, \Delta) \subset \Delta(\widetilde q, \Gamma).$
\end{proof}

In the rest of this section, we let $f := (f_1, \ldots, f_m) \colon \mathbb{R}^n \to \mathbb{R}^m$ be a polynomial mapping and $\Omega \subset \mathbb{R}^n$ an unbounded, closed and semialgebraic set. Let $n_i := \# (\Gamma(f_i) \cap \mathbb{Z}^n).$  For each vector $a := (a_1, \ldots, a_m) \in \mathbb{R}^{n_1} \times \cdots \times \mathbb{R}^{n_m},$ consider the polynomial mapping
$$F(a, \cdot) \colon \mathbb{R}^n \to \mathbb{R}^m, \quad x \mapsto (F_1(a_1, x), \ldots, F_m(a_m, x)),$$
where $F_i(a_i, x) := f_i(x) + \sum_{\kappa \in \Gamma(f_i)} a_{i, \kappa} x^\kappa$ for $i = 1, \ldots, m.$ 
With the above preparations, we have the following result, which establishes the openness of the nondegeneracy condition at infinity.

\begin{proposition}
Assume that the restriction $f|_\Omega$ is nondegenerate at infinity. Then for any vector $a \in \mathbb{R}^{n_1} \times \cdots \times \mathbb{R}^{n_m}$ with $\|a\|$ sufficiently small the restriction $F(a, \cdot)|_\Omega$ is nondegenerate at infinity.
\end{proposition}

\begin{proof}
Indeed, if the claim were not true, there would exist by Lemma~\ref{CurveSelectionLemma1}
a vector $q^0 \in \mathbb{R}^n$ with $\min_{j \in J} q^0_j < 0,$ nonempty sets $I_1 \subset \{1, \ldots, m\},$ $J \subset \{1, \ldots, n\},$ and analytic semialgebraic curves $c(t) \in (\mathbb{R}^*_+)^{I_1},$ $\eta(t) \in (\mathbb{R}^*)^J,$ and $a_i(t) \in \mathbb{R}^{n_i}$
with $\lim_{t \to 0} a_i(t) = 0$ for $i \in I_1,$ defined on a small enough interval $(0, \epsilon) \subset \mathbb{R},$ such that
\begin{eqnarray}
0 &>& d(q^0, \Gamma(f_i))  \ \textrm{ and } \ \Delta_i := \Delta(q^0, \Gamma(f_i)) \subset \mathbb{R}^J \quad \textrm{ for } i \in I_1, \nonumber\\
0 &=& F_{i, \Delta_i} (a_i(t), \eta(t)) \quad \textrm{ for } i \in I_1, \label{Eqn2}\\
0 &\in& \sum_{i \in I_1} c_i(t) \eta(t) \odot \nabla F_{i, \Delta_i}(a_i(t), \eta(t)) + C(\infty; \Omega) \cap \mathbb{R}^J. \label{Eqn3}
\end{eqnarray}

Let us expand $c_i(t), i \in I_1,$ and $\eta_j(t), j \in J,$ in terms of the parameter, say
\begin{eqnarray*}
c_i(t) &=& c_i^0 t^{k_i} + o(t^{k_i}) \quad \textrm{ and } \quad\eta_j(t) \ = \ x_j^0 t^{q_j} + o(t^{q_j}),
\end{eqnarray*}
where $c_i^0 > 0,$ $k_i \in \mathbb{Q},$ $x_j^0 \ne 0$ and $q_j \in {\mathbb Q}.$
Let $q := (q_1,\ldots,q_n)$ where $q_j := M$ for $j \not \in J$ with $M$ being sufficiently large and satisfying
\begin{eqnarray*}
M &>& \max_{i = 1, \ldots, m}\left\{\sum_{j \in J} q_j \kappa_j \mid \kappa \in {\Gamma}(f_i) \right\}.
\end{eqnarray*}
For each $i \in I_1,$ let
\begin{eqnarray*}
\widetilde{\Delta}_i & := & \Delta(q, \Delta_i) \ \subset \ \mathbb{R}^J.
\end{eqnarray*}
We first show that there exists a vector $\widetilde{q} \in \mathbb{R}^n$ with $\min_{j \in J} \widetilde{q}_j < 0$ such that $\widetilde{\Delta}_i = \Delta(\widetilde{q}, \Gamma(f_i)).$ In fact, if $q_j = 0$ for all $j \in J$ then $\widetilde{\Delta}_i = {\Delta}_i,$ and so the vector $\widetilde{q} := q^0$ has the desired properties.

Assume that $q_j \ne 0$ for some $j \in J.$ Then for $\delta > 0$ small, the vector $\widetilde q := q^0 + \delta \frac{q}{\|q\|}$ satisfies the condition
$d(\widetilde{q}, \Gamma(f_i))  < 0$ for all $i \in I_1.$ Moreover, in view of Lemma~\ref{PolyhedralLemma} we have for all $i \in I_1$ that
\begin{eqnarray*}
\Delta(\widetilde q, \Gamma(f_i)) &=& \Delta(\widetilde q,\Delta_i).
\end{eqnarray*}
On the other hand, we have for any $\kappa\in\Delta_i,$
\begin{eqnarray*}
\langle \widetilde q, \kappa\rangle &=& \langle q^0, \kappa \rangle + \frac{\delta}{\|q\|}\langle q,\kappa\rangle
\ = \ \min_{\kappa' \in \Gamma(f_i)}\langle q^0, \kappa'\rangle + \frac{\delta}{\|q\|}\langle q,\kappa\rangle
\ \geq \ \min_{\kappa' \in \Gamma(f_i)}\langle q^0, \kappa'\rangle + \frac{\delta}{\|q\|} \widetilde{d}_i,
\end{eqnarray*}
where $\widetilde{d}_i := \min_{\kappa' \in \Delta_i}\langle q, \kappa'\rangle;$ note that the equality happens if and only if $\kappa\in\widetilde\Delta_i,$ which turns out that $\widetilde\Delta_i = \Delta(\widetilde q,\Delta_i).$ Therefore, $\widetilde\Delta_i=\Delta(\widetilde q, \Gamma(f_i)),$ as required.

Recall that $a_i(t) \to 0$ as $t \to 0.$ Now a simple calculation shows that
\begin{eqnarray*}
F_{i, \Delta_i} (a_i(t), \eta(t)) &=& f_{i, \widetilde{\Delta}_i}(x^0) t^{\widetilde{d}_i} + o(t^{\widetilde{d}_i}),\\
\sum_{i \in I_1} c_i(t) \eta_j(t) \frac{\partial F_{i, \Delta_i}}{\partial x_j}(a_i(t), \eta(t))
&=& \sum_{i \in I_1} c^0_i x_j^0 \frac{\partial f_{i, \widetilde{\Delta}_i}}{\partial x_j}(x^0) t^{\widetilde{d}_i + {k_i}} + o(t^{\widetilde{d}_i + k_i}) \\
&=& \left( \sum_{i \in I_2} c^0_i x_j^0 \frac{\partial f_{i, \widetilde{\Delta}_i}}{\partial x_j}(x^0) \right) t^{\frak{m}} + o(t^{\frak{m}}),
\end{eqnarray*}
where $x^0 := (x_1^0, \ldots, x_n^0) \in ({\mathbb R}^*)^{J},$ with $x_j^0 = 0$ for $j \not \in J,$
$\frak{m} := \min_{i \in I_1} (\widetilde{d}_i + {k_i}),$ and $I_2 := \{i \in I_1 \mid \widetilde{d}_i + {k_i} = \frak{m}\}.$
These equalities, in combining with \eqref{Eqn2} and \eqref{Eqn3}, tell us that
\begin{eqnarray*}
 0 &=& f_{i, \widetilde{\Delta}_i}(x^0) \quad \textrm{ for all } \quad i \in I_1, \\
0 &\in& \sum_{i \in I_2} c^0_i x^0 \odot \nabla f_{i, \widetilde{\Delta}_i}(x^0) + o(1) + C(\infty; \Omega) \cap \mathbb{R}^J.
\end{eqnarray*}
Letting $t \to 0$ and noting that $C(\infty; \Omega)$ is a closed cone, the second relation yields
\begin{eqnarray*}
0 &\in&
\sum_{i \in I_2} c^0_i x^0 \odot \nabla f_{i, \widetilde{\Delta}_i}(x^0) + C(\infty; \Omega) \cap \mathbb{R}^J,
\end{eqnarray*}
which contradicts the nondegeneracy at infinity of $f|_\Omega.$
\end{proof}

\section{Continuous selections defined by nondegenerate mappings at infinity} \label{Section4}

In this section, let $f := (f_1, \ldots, f_m) \colon \mathbb{R}^n \to \mathbb{R}^m$ be a polynomial mapping and $\Omega \subset \mathbb{R}^n$ an unbounded, closed and semialgebraic set. Following Kouchnirenko \cite{Kouchnirenko1976}, we say that $f$ is {\em convenient} if its component polynomials of $f$ are convenient, i.e.,  for every $i = 1, \ldots, m,$ the Newton polyhedron $\Gamma(f_i)$ of $f_i$ intersects each coordinate axis at a point different from the origin $0$ in $\mathbb{R}^n.$

The following result plays a crucial role in later developments and is inspired by the work of Broughton~\cite{Broughton1988}, where it was proven that a complex polynomial $P \colon \mathbb{C}^n \to \mathbb{C},$ which is convenient and nondegenerate at infinity, is tame in the sense that there exists a compact neighborhood $U$ of the critical points of $P$ such that $\|\nabla P(x)\|$ is bounded away from $0$ on $\mathbb{C}^n \setminus U.$

\begin{theorem}\label{Thm41}
Assume that the mapping $f$ is convenient and the restriction $f|_\Omega$ is nondegenerate at infinity. Let $h \colon \mathbb{R}^n \to \mathbb{R}$ be a continuous selection of the polynomials $f_1, \ldots, f_m.$
Then the following statements hold$:$
\begin{enumerate}[{\rm (i)}]
\item there exists a constant $\delta > 0$ such that for all $x \in \Omega$ with $\|x\|$ sufficiently large, all $u \in \partial h(x)$ and all $v \in N(x; \Omega),$ we have
\begin{eqnarray*}
\|u + v \| &\ge& \delta;
\end{eqnarray*}

\item if $h$ is bounded from below on $\Omega,$ then there exists a constant $c > 0$ such that for all $x \in \Omega$ with $\|x\|$ sufficiently large, we have
\begin{eqnarray*}
h(x) &\ge& c\, \|x\|,
\end{eqnarray*}
in particular, $h$ is coercive on $\Omega$ in the sense that
\begin{eqnarray*}
\liminf_{\|x\| \to \infty, x \in \Omega} h(x) &=& \infty.
\end{eqnarray*}
\end{enumerate}
\end{theorem}

\begin{proof}
We first note that the function $h$ is locally Lipschitz (by \cite[Proposition~4.1.2]{Scholtes2012}) and semialgebraic (by \cite[Theorem~3.1]{Guo2024}).

(i) Assume for contradiction that there exist sequences $x_k \in \Omega,$ $u_k \in \partial h(x_k)$ and $v_k \in N(x_k; \Omega)$ such that
$\|x_k\| \to \infty$ and $\|u_k + v_k\| \to 0$ as $k \to \infty.$

In view of Lemma~\ref{Lemma24}, $\partial h(x_k)$ is a subset of the convex hull in $\mathbb{R}^n$ of the vectors $\nabla f_i(x_k)$ with $f_i(x_k) =h(x_k).$ Hence, passing to subsequences if necessary, we can find a nonempty set $I_1 \subset \{1, \ldots, m\}$ and a vector $\lambda_k  \in \mathbb{R}^m$ satisfying the following conditions
\begin{eqnarray*}
&& u_k = \sum_{i \in I_1} \lambda_{i, k} \, \nabla f_i(x_k), \quad \sum_{i \in I_1} \lambda_{i, k} = 1,\\
&& \lambda_{i, k}  > 0 \quad \textrm{ and } \quad f_i(x_k)\ = \ h(x_k) \ \ \textrm{ for all } i \in I_1.
\end{eqnarray*}
(Of course, it can happen that $f_i(x_k) =  h(x_k)$ for some $i \not \in I_1.$)

By Lemma~\ref{CurveSelectionLemmaAtInfinity}, there exists an analytic semialgebraic curve
$$(\eta, v, \lambda) \colon (0, \epsilon) \rightarrow \mathbb{R}^n \times \mathbb{R}^n \times \mathbb{R}^m, \quad t \mapsto \big (\eta(t), v(t), \lambda(t)\big),$$
satisfying the following conditions
\begin{enumerate}[\upshape (c1)]
\item $\eta(t) \in \Omega$ and $v(t) \in N(\eta(t); \Omega);$
\item $\lim_{t \to 0^+} \|\eta(t)\|= \infty;$
\item $f_i(\eta(t)) \equiv h(\eta(t))$ for all $i \in I_1;$
\item $u(t) := \sum_{i \in I_1} \lambda_i(t) \nabla f_i(\eta(t)) \in \partial h(\eta(t));$
\item $\lim_{t \to 0^+}\|u(t) + v(t) \| = 0;$
\item $\sum_{i \in I_1} \lambda_i(t) \equiv 1$ and $\lambda_i(t) > 0$ for all $i \in I_1.$
\end{enumerate}

Let $J := \{j \mid \eta_j \not\equiv 0\}.$ By Condition~(c2), the set $J$ is nonempty. In view of Lemma~\ref{GrowthDichotomyLemma}, for each $j \in J,$ the coordinate $\eta_j$ is expanded as follows
\begin{eqnarray*}
\eta_j(t) &=& {x}^0_j \, t^{q_j} + \cdots,
\end{eqnarray*}
where ${x}^0_j \ne 0$ and $q_j \in \mathbb{Q}.$ Here and in the following, the dots stand for the higher order terms in $t.$
From Condition~(c2), we get $q_* := \min_{j \in J} q_j< 0.$

Let $q_j := M$ for $j \not \in J$ with $M$ being sufficiently large and satisfying
\begin{eqnarray*}
M &>& \max_{i = 1, \ldots, m}\left\{\sum_{j' \in J} q_{j'} \kappa_{j'} \mid \kappa \in {\Gamma}(f_i) \right\}.
\end{eqnarray*}
For each $i = 1, \ldots, m,$ let $d_i$ be the minimal value of the linear function $\langle q, \kappa\rangle$ on ${\Gamma}(f_i)$ and let $\Delta_i$ be the maximal face of ${\Gamma}(f_i)$ (maximal with respect to the inclusion of faces) where the linear function takes this value, i.e.,
\begin{equation*}
d_i \ := \ d(q, {\Gamma}(f_i)) \quad \textrm{ and } \quad \Delta_i \ := \ \Delta(q, {\Gamma}(f_i)).
\end{equation*}
Since $f_i$ is convenient, for each $j=1,\ldots,n$, there exists an integer $\ell_{ij}\geq 1$ such that $\ell_{ij}e_j\in\Gamma(f_i)$, where $e_j$ denotes the $j$th canonical basis vector of $\mathbb{R}^n.$
Recall that
\begin{eqnarray*}
q_* &:=& \min_{j\in J}q_j \ < \ 0.
\end{eqnarray*}
Choose $j_*\in J$ such that $q_{j_*} = q_*.$  Then, for all $j = 1, \ldots, n,$
\begin{eqnarray}\label{Eqn4}
d_i &=& d(q,\Gamma(f_i)) \ \leq \ \langle q,\ell_{ij_*}e_{j_*} \rangle \ = \ \ell_{ij_*}q_* \ \leq \ q_* \ \leq \ q_j.
\end{eqnarray}

Recall that $\mathbb{R}^J := \{x := (x_1, \ldots, x_n) \in \mathbb{R}^n \mid x_j = 0 \textrm{ for } j \not \in J \}.$
Take any $i \in I_1.$ Since $f_i$ is convenient, the restriction of $f_i$ to $\mathbb{R}^J$ is not constant and ${\Gamma}(f_i) \cap {\mathbb{R}^J} = {\Gamma}(f_i|_{\mathbb{R}^J}) $ is nonempty and different from $\{0\}.$ Furthermore, we have from the definition of the vector $q$ that
\begin{equation*}
d_i = d(q, {\Gamma}(f_i|_{\mathbb{R}^J})) \quad \textrm{ and } \quad \Delta_i = \Delta(q, {\Gamma}(f_i|_{\mathbb{R}^J}))
\ \subset \ {\mathbb{R}^J}.
\end{equation*}
Consequently, for each $j \not \in J,$ the polynomial $f_{i, \Delta_i}$ does not depend on the variable $x_j.$
Now suppose that $f_i$ is written as $f_i(x) = \sum_{\kappa} a_{i, \kappa} x^\kappa.$ A simple calculation shows that
\begin{eqnarray*}
f_i(\eta(t)) &=& \sum_{\kappa\in{\Gamma}(f_i) \cap \mathbb{R}^J} a_{i, \kappa} (\eta(t))^\kappa\\
&=& \sum_{\kappa \in {\Gamma}(f_i) \cap \mathbb{R}^J} \Big ( a_{i, \kappa} \prod_{j = 1}^n \eta_j(t)^{\kappa_j} \Big) \\
&=& \sum_{\kappa \in {\Gamma}(f_i) \cap \mathbb{R}^J}\Big (a_{i, \kappa} \prod_{j \in J} ({x}^0_j t^{q_j})^{\kappa_j}+ \cdots \Big)\\
&=& \sum_{\kappa \in {\Gamma}(f_i) \cap \mathbb{R}^J}\Big(a_{i, \kappa} ({x}^*)^\kappa t^{\sum_{j \in J}q_j \kappa_j} +\cdots \Big)\\
&=& \sum_{\kappa \in \Delta_i}a_{i, \kappa} ({x}^*)^\kappa t^{d_i} + \cdots,
\end{eqnarray*}
where $x^* := (x^*_1, \ldots, x^*_n) \in (\mathbb{R}^*)^J$ is defined by
\begin{eqnarray*}
x^*_j &:=&
\begin{cases}
0 & \textrm{ if } j \not \in J, \\
x^0_j& \textrm{ otherwise.}
\end{cases}
\end{eqnarray*}
Recall that $f_{i, \Delta_i} (x) = \sum_{\kappa \in \Delta_i} a_{i, \kappa} x^\kappa.$ Hence
\begin{eqnarray*}
f_i(\eta(t)) &=& f_{i, \Delta_i}({x}^*) \, t^{d_i} + \cdots.
\end{eqnarray*}

Take any $j \in J.$ We have
\begin{eqnarray*}
\frac{\partial f_i}{\partial x_j}(\eta(t))
&=& \frac{\partial f_{i, \Delta_i}}{\partial x_j}({x}^*) \, t^{d_i - q_j}+ \cdots.
\end{eqnarray*}
On the other hand, from Condition~(c6) and Lemma~\ref{GrowthDichotomyLemma} we can write
\begin{eqnarray*}
\lambda_i(t) &=& \lambda^0_i \, t^{k_i} + \cdots
\end{eqnarray*}
for some $\lambda^0_i > 0$ and ${k_i} \ge 0.$ Therefore,
\begin{eqnarray*}
u_j(t)
&=& \sum_{i \in I_1} \lambda_i^0 \frac{\partial f_{i, \Delta_i}}{\partial x_j}({x}^*) \, t^{{k_i} + d_i - q_j}+ \cdots\ = \
\left( \sum_{i \in I_2} \lambda_i^0 \frac{\partial f_{i, \Delta_i}}{\partial x_j}({x}^*) \right) t^{\mathfrak{m} - q_j}+ \cdots,
\end{eqnarray*}
where $\mathfrak{m}:= \min_{i \in I_1} ({k_i} + d_i)$ and $I_2 := \{i \in I_1 \mid {k_i} + d_i = \mathfrak{m}\} \ne \emptyset.$
It follows from Condition~(c6) that $k_{i_*} = 0$ for some index $i_* \in {I}_1.$
In combining this with \eqref{Eqn4}, we obtain $\mathfrak{m} \le d_{i_*} \le q_j.$

If $v_j(t) \ne 0$ for all $t > 0$ sufficiently small, then it follows from Lemma~\ref{GrowthDichotomyLemma} that
\begin{eqnarray*}
v_j(t) &=&{v}^0_j \, t^{p_j} + \cdots
\end{eqnarray*}
for some $v^0_j \ne 0$ and $p_j \in \mathbb{Q},$ otherwise we put $v^0_j= 0$ and $p_j := \infty.$

From Condition~(c5), we know that $u_j(t) + v_j(t) \to 0$ as $t \to 0^+.$ Therefore,
\begin{eqnarray*}
\sum_{i \in I_2} \lambda_i^0 \frac{\partial f_{i, \Delta_i}}{\partial x_j}({x}^*) &=&
\begin{cases}
0 & \textrm{ if } p_j> \mathfrak{m} - q_j, \\
-v^0_j& \textrm{ otherwise.}
\end{cases}
\end{eqnarray*}

We next put $\xi := (\xi_1, \ldots, \xi_n) \in \mathbb{R}^J,$ where
\begin{eqnarray*}
\xi_j &:=&
\begin{cases}
0 & \textrm{ if } \ j \not \in J \ \textrm{ or } \ p_j> \mathfrak{m} - q_j, \\
x^0_j v^0_j& \textrm{ otherwise.}
\end{cases}
\end{eqnarray*}
Certainly
\begin{eqnarray} \label{Eqn5}
\sum_{i \in I_2} \lambda^0_i \,x^* \odot\nabla f_{i, \Delta_i} (x^*) + \xi &= & 0.
\end{eqnarray}
Observe that if $\xi \ne 0,$ then $\eta(t) \odot v(t) \ne 0$ for all $t > 0$ small enough, and a direct computation shows that
\begin{eqnarray*}
\lim_{t \to 0^+} \frac{\eta(t) \odot v(t)}{\|\eta(t) \odot v(t)\|} &=& \frac{\xi}{\|\xi\|}.
\end{eqnarray*}
Then we conclude from Lemma~\ref{MD32} and Condition~(c1) that $\xi \in C(\infty; \Omega).$

Since the restriction $f|_\Omega$ is nondegenerate at infinity, it follows from \eqref{Eqn5} that $f_{i, \Delta_i}(x^*) \ne 0$ for some $i \in I_2.$ Recall from Condition~(c3) that
\begin{eqnarray*}
h(\eta(t)) &=& f_i(\eta(t)) \ = \ f_{i, \Delta_i}({x}^*) \, t^{d_i} + \cdots,
\end{eqnarray*}
which yields
\begin{eqnarray*}
\frac{d}{dt} (h \circ \eta) (t) &=&d_i f_{i, \Delta_i}({x}^*) \, t^{d_i - 1} + \cdots.
\end{eqnarray*}

On the other hand, we know from Condition~(c1) that
$(\delta_\Omega \circ \eta)(t) = 0$ and $v(t) \in N(\eta(t); \Omega) = \partial \delta_\Omega (\eta(t))$ (see Lemma~\ref{Lemma23}). By Lemma~\ref{ChainRule}, we have for all $t > 0$ sufficiently small,
\begin{eqnarray*}
\frac{d}{dt} (h \circ \eta) (t) &=&\frac{d}{dt} (h \circ \eta) (t)+ \frac{d}{dt} (\delta_\Omega \circ \eta) (t) \\
&=&\left \langle u(t), \frac{d \eta(t)}{dt} \right \rangle + \left \langle v(t), \frac{d \eta(t)}{dt} \right \rangle \\
&=&\left \langle u(t) + v(t), \frac{d \eta(t)}{dt} \right \rangle.
\end{eqnarray*}
Consequently,
\begin{eqnarray*}
\left| \frac{d}{dt} (h \circ \eta) (t)\right|&\le& \|u(t) + v(t) \| \left\| \frac{d \eta(t)}{dt} \right\|.
\end{eqnarray*}
Therefore,
\begin{eqnarray*}
c\, t^{d_i - q_*} &\le& \|u(t) + v(t) \|
\end{eqnarray*}
for some $c > 0,$ which contradicts Condition~(c5) because $d_i - q_* \le 0,$ as seen in \eqref{Eqn4}. The proof of the statement (i) is complete.

(ii) We first observe from Corollary~\ref{Corollary214} that for all $s$ sufficiently large, $\Omega \cap \mathbb{S}_s$ is a nonempty compact set. By the Weierstrass theorem and Lemma~\ref{CurveSelectionLemmaAtInfinity}, we can find an analytic semialgebraic curve $\eta \colon (R, \infty) \to \mathbb{R}^n$ with $\eta(s) \in \Omega \cap \mathbb{S}_s$ such that $\eta(s)$ is an optimal solution of the problem
\begin{eqnarray*}
\textrm{minimize} \ h(x) \quad \textrm{ subject to } \quad {x \in \Omega \cap \mathbb{S}_s}.
\end{eqnarray*}
The function $h$ is locally Lipschitz. By Lemma~\ref{FermatRule}, hence
\begin{eqnarray*}
0 &\in& \partial h(\eta(s)) + N(\eta(s); \Omega \cap \mathbb{S}_s).
\end{eqnarray*}
Observe from Corollary~\ref{Corollary214} that for all $s$ sufficiently large,
\begin{eqnarray*}
N(\eta(s); \Omega \cap \mathbb{S}_s) &\subset& N(\eta(s); \Omega) + N(\eta(s); \mathbb{S}_s).
\end{eqnarray*}
Therefore, there exist $u(s) \in \partial h(\eta(s)),$ $v(s) \in N(\eta(s); \Omega)$ and $\mu(s) \in \mathbb{R}$ such that
\begin{eqnarray*}
0 &=& u(s) + v(s) + \mu(s) \eta(s) ,
\end{eqnarray*}
which, together with (i), yields
\begin{eqnarray*}
|\mu(s) s| &=& |\mu(s)| \|\eta(s)\| \ = \ \|u(s) + v(s)\| \ \ge \ \delta.
\end{eqnarray*}
On the other hand, it follows from Lemma~\ref{ChainRule} that for all $s$ sufficiently large,
\begin{eqnarray*}
\frac{d}{ds} (h \circ \eta) (s) &=&\frac{d}{ds} (h \circ \eta) (s)+ \frac{d}{ds} (\delta_\Omega \circ \eta) (s) \\
&=&\left \langle u(s), \frac{d \eta(s)}{ds} \right \rangle + \left \langle v(s), \frac{d \eta(s)}{ds} \right \rangle \\
&=&\left \langle u(s) + v(s), \frac{d \eta(s)}{ds} \right \rangle \\
&=& - \mu(s) \left \langle \eta(s), \frac{d \eta(s)}{ds} \right \rangle \\
&=& - \mu(s) \frac{1}{2} \frac{d}{ds} \|\eta(s)\|^. \\
&=& - \mu(s) s.
\end{eqnarray*}
The function $h \circ \eta$ is semialgebraic, and so it is either decreasing or increasing in view of Lemma~\ref{MonotonicityLemma} (perhaps after increasing $R$). In the first case, we have
\begin{eqnarray*}
\frac{d}{ds} (h \circ \eta) (s) &=& - \left| \frac{d}{ds} (h \circ \eta) (s) \right| \ = \ - |\mu(s) s| \ \le \ -\delta.
\end{eqnarray*}
which yields that $(h \circ \eta) (s) \to -\infty$ when $s \to \infty,$ in contradiction to our assumption that $h$ is bounded from below. Therefore,
the function $h \circ \eta$ is increasing and satisfies
\begin{eqnarray*}
\frac{d}{ds} (h \circ \eta) (s) &=& \left| \frac{d}{ds} (h \circ \eta) (s) \right| \ = \ |\mu(s) s| \ \ge \ \delta.
\end{eqnarray*}
By letting $c := \frac{\delta}{2} > 0,$ we get $h(\eta(s)) \ge c s$ for all $s$ sufficiently large.
Finally, for any $x \in \Omega$ with $s := \|x\|$ sufficiently large, we have
\begin{eqnarray*}
h(x) &\ge& h(\eta(s)) \ \ge \ c s \ = \ c\, \|x\|,
\end{eqnarray*}
as required. The proof is complete.
\end{proof}



Theorem~\ref{Thm41} allows us to prove a nonconvex version of \cite[Theorem~1]{Luo1999} and \cite[Corollary~1]{Belousov2002}. Define for $\varepsilon := (\epsilon_1, \ldots, \epsilon_m) \in \mathbb{R}^m,$
\begin{eqnarray*}
X(\varepsilon) &:=& \{x \in \Omega \mid f_i(x) \le \epsilon_i, \ i = 1, \ldots, m\}.
\end{eqnarray*}

\begin{corollary}
Assume that the mapping $f$ is convenient and the restriction $f|_\Omega$ is nondegenerate at infinity.
If the sets $X(\varepsilon_k)$ are nonempty for some sequence of nonnegative vectors $\varepsilon_k$ approaching zero, then the set $X(0)$ is also nonempty.
\end{corollary}
\begin{proof}
Assume for contradiction that $X(0) = \emptyset.$ Then the function $h(x) := \max_{i = 1, \ldots, m} f_i(x)$ is positive on $\Omega;$ moreover, since
the sets $X(\varepsilon_k)$ are nonempty,
 there exists a sequence $x_k \in \Omega$ such that $\|x_k\| \to \infty$ and $h(x_k) \to 0.$ On the other hand, by definition, $h$ is a continuous selection of the polynomials $f_1, \ldots, f_m.$ Putting these facts together, we get a contradiction to Theorem~\ref{Thm41}(ii).
\end{proof}

\section{Pareto optimality under nondegeneracy at infinity} \label{Section5}

In the rest of the paper, we consider the multiobjective optimization problem formulated in Section~\ref{Section1}:
\begin{equation}\label{MP}
{\rm Min}_{\mathbb{R}_+^m} \left\{ f(x) \mid x \in \Omega\right\} \tag{MP}
\end{equation}
with $f := (f_1, \ldots, f_m) \colon \mathbb{R}^n \to \mathbb{R}^m$ being a polynomial mapping and $\Omega \subset \mathbb{R}^n$ being an 
unbounded, closed and semialgebraic set. 

\subsection{Some preparations}

In the scalar optimization case (i.e., $m = 1$) it is easy to see that the following conditions are equivalent:
\begin{enumerate}[\upshape (i)]
\item $f$ is bounded from below on $\Omega;$
\item there exists $\bar z \in \mathbb{R}$ such that $f(\Omega) \cap (\bar z - \mathbb{R}_+)$ is an empty set.
\item there exists $\bar z \in \mathbb{R}$ such that $f(\Omega) \cap (\bar z - \mathbb{R}_+)$
is a nonempty and bounded set.
\end{enumerate}
These equivalences lead to the following definition in the multiobjective optimization case (i.e., $m \geq 1).$

\begin{definition}{\rm
\begin{enumerate}[\upshape (i)]
\item The restriction $f|_\Omega$ is said to have {\em an empty section} if there exists $\bar z \in \mathbb{R}^m$ such that
$f(\Omega) \cap \left(\bar z - \mathbb{R}_+^m\right)$ is an empty set.

\item The restriction $f|_\Omega$ is said to have  {\em a nonempty bounded section} if there exists $\bar z \in \mathbb{R}^m$ such that
$f(\Omega) \cap \left(\bar z - \mathbb{R}_+^m\right)$ is a nonempty and bounded set.
\end{enumerate}
}\end{definition}

By definition, if the restriction $f|_{\Omega}$ has a nonempty bounded section then it has an empty section, but not conversely (see Example~\ref{Example67}). Moreover, it is not hard to check that the condition that the restriction $f|_{\Omega}$ has an empty section (resp., a nonempty bounded section) is necessary (but not sufficient) for the problem~\eqref{MP} to admit a weakly efficient solution (resp., an efficient solution).
Interestingly, these conditions are also sufficient for a large class of optimization problems with semialgebraic data; see Theorem~\ref{Thm52} below.

Recall from \cite{Wierzbicki1986} (see also \cite[Section 4.6]{Ehrgott2005} and \cite{PHAMTS2018-1}) that the {\em achievement function} $\psi\colon \mathbb{R}^n \to \mathbb{R} \cup \{\infty\}$ for the problem~\eqref{MP} is defined by
\begin{eqnarray*}
\psi(x) &:=& \sup_{y \in \Omega} \min_{i=1,\ldots,m} [f_i(x) - f_i(y)].
\end{eqnarray*}
Then the function $\psi$ is nonnegative on $\Omega$ (by definition), semialgebraic (by the Tarski--Seidenberg theorem), and moreover, it completely determines the set of all weakly efficient solutions to the problem~\eqref{MP}, denoted by $\mathrm{WSOL},$ in the sense that
\begin{eqnarray*}
\mathrm{WSOL} &=& \{ x \in \Omega \mid \psi(x) = 0\}.
\end{eqnarray*}
Let us now provide necessary and sufficient conditions for the restriction $f|_\Omega$ to have an empty section, which is of its own interest.

\begin{theorem}\label{Thm51}
The following statements are equivalent$:$
\begin{enumerate}[\upshape (i)]
\item the restriction $f|_\Omega$ has an empty section$;$
\item the function $x \mapsto \max_{i = 1, \ldots, m} f_i(x)$ is bounded from below on $\Omega;$
\item $\psi(x) < \infty$ for all $x \in \mathbb{R}^n;$
\item $\psi(\bar x) < \infty$ for some $\bar x \in \mathbb{R}^n;$
\item $\inf_{x \in \Omega} \psi(x) = 0.$
\end{enumerate}
\end{theorem}
\begin{proof}
(i) $\Rightarrow$ (ii)
Let $\bar z \in \mathbb{R}^m$ be such that $f(\Omega) \cap \left(\bar z - \mathbb{R}_+^m\right) = \emptyset,$ and let $x \in \Omega$ be arbitrary.
Then $f(x) \not\in \bar z - \mathbb{R}_+^m,$ and so $f_i(x) \geq \bar z_i$ for some index $i.$
Hence
\begin{eqnarray*}
\max_{i = 1, \ldots, m}f_i(x) &\geq& \min_{i = 1, \ldots, m}\bar z_i.
\end{eqnarray*}
Since this is true for arbitrary $x \in \Omega,$ the required statement follows.

(ii) $\Rightarrow$ (iii) Assume that the function $\max_{i = 1, \ldots, m}f_i$ is bounded from below on $\Omega,$ and let $x \in \mathbb{R}^n$ be arbitrary.
We have for all $y \in \Omega,$
\begin{eqnarray*}
\max_{i = 1, \ldots, m}[f_i(y) - f_i(x)] &\geq& \max_{i = 1, \ldots, m}f_i(y) - \max_{i = 1, \ldots, m} f_i(x).
\end{eqnarray*}
Consequently,
\begin{eqnarray*}
- \psi(x) &=& \inf_{y \in \Omega}\max_{i = 1, \ldots, m}[f_i(y) - f_i(x)] \\
&\geq& \inf_{y \in \Omega}\max_{i = 1, \ldots, m}f_i(y) - \max_{i = 1, \ldots, m} f_i(x)\\
&>& -\infty,
\end{eqnarray*}
which yields $\psi(x) < \infty.$
Since this is true for arbitrary $x \in \mathbb{R}^n,$ the required statement follows.

(iii) $\Rightarrow$ (iv) This implication is trivial.

(iv) $\Rightarrow$ (v)
Assume that $\psi (\bar x ) < \infty$ for some $\bar x \in \mathbb{R}^n.$
By definition, $\psi(x) \geq 0$ for all $x \in \Omega,$ and so $\delta := \inf_{x \in \Omega} \psi(x) \geq 0.$
Assume for contradiction that $\delta > 0.$ Then for each $x \in \Omega,$ we have $\psi(x) \geq \delta,$ and so there exists $y_x \in \Omega$ such that $\min_{i=1,\ldots,m}[f_i(x) - f_i(y_x)] > \tfrac{\delta}{2}.$

Now, fix any $x_0 \in \Omega.$ There exists $x_1 := y_{x_0} \in \Omega$ such that
\begin{eqnarray*}
\min_{i=1,\ldots,m}[f_i(x_0) - f_i(x_1)] &>& \tfrac{\delta}{2}.
\end{eqnarray*}
Similarly, for $x_1$, there exists $x_2 := y_{x_1} \in \Omega$ such that
\begin{eqnarray*}
\min_{i=1,\ldots,m}[f_i(x_1) - f_i(x_2)] &>& \tfrac{\delta}{2}.
\end{eqnarray*}
By induction, we construct a sequence $x_k \in \Omega$ such that for all $k \geq 1,$
\begin{eqnarray*}
\min_{i=1,\ldots,m}[f_i(x_{k-1}) - f_i(x_k)] &>& \tfrac{\delta}{2}.
\end{eqnarray*}
Adding the above inequalities, we get
\begin{eqnarray*}
\min_{i=1,\ldots,m}[f_i(x_0) - f_i(x_k)] &>& k \cdot \tfrac{\delta}{2}.
\end{eqnarray*}
It follows that
\begin{eqnarray*}
\min_{i = 1, \ldots, m}[f_i(\bar{x}) - f_i(x_k)] &\geq& \min_{i=1,\ldots,m}[f_i(\bar{x})-f_i(x_0)]+\min_{i=1,\ldots,m}[f_i(x_0)-f_i(x_k)] \\
&>& \min_{i=1,\ldots,m}[f_i(\bar{x})-f_i(x_0)]+\frac{k\delta}{2}.
\end{eqnarray*}
By definition, then
\begin{eqnarray*}
\psi(\bar x) &>& \min_{i=1,\ldots,m}[f_i(\bar{x})-f_i(x_0)]+\frac{k\delta}{2}.
\end{eqnarray*}
Letting $k \to \infty,$ we obtain $\psi(\bar x) = \infty,$ a contradiction.

(v) $\Rightarrow$ (i)
Since $\inf_{x \in \Omega} \psi(x) = 0,$ we have $\psi(\bar x) < \infty$ for some $\bar x \in \mathbb{R}^n.$
Now, define $\bar z := \left(\bar z_1, \ldots, \bar z_m\right) \in \mathbb{R}^m$ by
\begin{eqnarray*}
\bar z &:=& f(\bar x) - \psi(\bar x) (1, \ldots, 1).
\end{eqnarray*}
For arbitrary $y \in \Omega,$ we have by definition that
\begin{eqnarray*}
\max_{i = 1, \ldots, m}[f_i(y) - f_i(\bar x)] & \geq & -\psi(\bar x).
\end{eqnarray*}
Then, there exists some $i$ such that
\begin{eqnarray*}
f_i(y) - f_i(\bar x) &\geq& -\psi(\bar x);
\end{eqnarray*}
or equivalently,
\begin{eqnarray*}
f_i(y) &\geq& f_i(\bar x)- \psi(\bar x) \ = \ \bar z_i.
\end{eqnarray*}
Consequently, $f(y) \not\in \bar z - \mathbb{R}^m_+.$
Since this is true for arbitrary $y \in \Omega,$ we get
\begin{eqnarray*}
f(\Omega) \cap \left(\bar z - \mathbb{R}_+^m\right)& = & \emptyset,
\end{eqnarray*}
as desired.
\end{proof}

\begin{remark}{\rm
(i) By the Tarski--Seidenberg theorem, the function $\max_{i = 1, \ldots, m} f_i$ is semialgebraic. On the other hand, we observe that the recent results in \cite{Lee2026, PHAMTS2023} can assist in detecting whether an extended-real-valued semialgebraic function is bounded from below. These facts, together with Theorem~\ref{Thm51}, show that the condition that the restriction $f|_\Omega$ has an empty section can be checked.

(ii) We would like to note that necessary and sufficient conditions for the existence of nonempty bounded sections are not discussed yet.
This interesting question will be investigated in the future work.
}\end{remark}

Now we consider the function $\phi \colon \mathbb{R}^n \times \mathbb{R}^n \to \mathbb{R}$ by defining
\begin{eqnarray*}\label{Eqn6}
\phi(x, y) &:=& \max_{i = 1, \ldots, m} [f_i(y) - f_i(x)].
\end{eqnarray*}
For each $x \in \mathbb{R}^n,$ let
\begin{eqnarray*}
\Omega(x) &:=& \left\{y \in \Omega \mid \phi (x, y) = - \psi(x) \right\},
\end{eqnarray*}
which is just the set of minimizers of the optimization problem $\inf_{y \in \Omega} \phi(x, y)$ because we have by definition that
\begin{eqnarray*}
\inf_{y \in \Omega} \phi(x, y) &=& - \sup_{y \in \Omega} [-\phi(x, y)] \ = \ - \psi(x).
\end{eqnarray*}
Recall that the set of weakly efficient solutions to the problem~\eqref{MP} is expressed in terms of $\psi$ by
\begin{eqnarray*}
\mathrm{WSOL} &=& \left\{ x \in \Omega \mid \psi(x) = 0\right\}.
\end{eqnarray*}
The following lemma is well known (see, for example, \cite{Ehrgott2005, PHAMTS2018-1}); we include its proof for the sake of completeness.

\begin{lemma}\label{Lem51}
The following statements hold$:$
\begin{enumerate}[\upshape (i)]
\item the function $\psi$ is finite, locally Lipschitz and semialgebraic provided that it is finite at some point $\bar{x} \in \mathbb{R}^n;$

\item the function $\phi$ is locally Lipschitz, semialgebraic and satisfies the triangle inequality
\begin{eqnarray*}
\phi(x, y) &\le& \phi(x, z) + \phi(z, y)   \quad \textrm{ for all } \quad x, y, z \in \mathbb{R}^n;
\end{eqnarray*}

\item for any $x \in \mathbb{R}^n,$ $\Omega(x)$ is a (possibly empty) subset of $\mathrm{WSOL}.$
\end{enumerate}
\end{lemma}

\begin{proof}
(i) Assume that the function $\psi$ is finite at some point $\bar{x} \in \mathbb{R}^n.$ By Theorem~\ref{Thm51}, $\psi$ is finite on $\mathbb{R}^n,$ then it is locally Lipschitz (by a direct checking) and semialgebraic (by the Tarski--Seidenberg theorem).

(ii) The function $\phi$ is locally Lipschitz because it is the pointwise maximum of finitely many polynomial functions; $\phi$ is semialgebraic by the Tarski--Seidenberg theorem. To see the triangle inequality, take any $x, y, z \in \mathbb{R}^n.$ We have
\begin{eqnarray*}
\phi(x, y) &=& \max_{i = 1, \ldots, m} [f_i(y) - f_i(x)] \\
&=& \max_{i = 1, \ldots, m} [f_i(y) - f_i(z) + f_i(z) - f_i(x)] \\
&\leq& \max_{i = 1, \ldots, m}[f_i(y) - f_i(z)] + \max_{i = 1, \ldots, m} [f_i(z) - f_i(x)] \\
&=& \phi(z, y) + \phi(x, z).
\end{eqnarray*}

(iii) Fix any $x \in \mathbb{R}^n$ and let $y \in \Omega(x).$ We have for any $y' \in \Omega,$
\begin{eqnarray*}
\phi(x, y') - \phi(x, y) &\ge& 0,
\end{eqnarray*}
and hence by (ii), $\phi(y, y') \ge 0.$ This is true for any $y' \in \Omega,$ so it implies that
\begin{eqnarray*}
-\psi(y) &=& \inf_{y' \in \Omega} \phi(y, y') \ \ge \ 0,
\end{eqnarray*}
from which follows $\psi(y) = 0,$ or equivalently, $y \in \mathrm{WSOL},$ as required.
\end{proof}

\subsection{The existence of efficient solutions and weak sharp minima}

For simplicity in the following, we denote by $\mathrm{ESOL}$ and $\mathrm{EVAL}$ (resp., $\mathrm{WSOL}$ and $\mathrm{WVAL}$) the sets of efficient solutions and efficient values (resp., the set of weakly efficient solutions and weakly efficient values) to the problem ~\eqref{MP}. For each $z \in \mathbb{R}^m,$ we put ${L}({z}) := f^{-1}({z} - \mathbb{R}^m_+) \cap \Omega.$
\begin{theorem} \label{Thm52}
Assume that the mapping $f$ is convenient and the restriction $f|_\Omega$ is nondegenerate at infinity. Then the following conditions are equivalent$:$
\begin{enumerate}[{\rm (c1)}]
\item the set $\mathrm{ESOL}$ of all efficient solutions to the problem~\eqref{MP} is nonempty$;$

\item the set $\mathrm{WSOL}$ of all weakly efficient solutions to the problem~\eqref{MP} is nonempty$;$

\item the restriction $f|_\Omega$ has an empty section$;$

\item for any $\bar{x} \in \mathbb{R}^n,$ there exist constants $c_1 > 0$ and $c_2 \in \mathbb{R}$ such that
\begin{eqnarray*}
\phi({\bar{x}}, y) &\ge& c_1 \|y\| + c_2 \quad \textrm{ for all } \quad y \in \Omega;
\end{eqnarray*}

\item there exists $\bar{x} \in \mathbb{R}^n$ such that the function $\phi({\bar{x}}, \cdot)$ is coercive on $\Omega;$

\item for any $\bar{z} \in \mathbb{R}^m,$ the set $L(\bar{z})$ is compact$;$

\item the restriction $f|_\Omega$ has a nonempty bounded section.
\end{enumerate}
When these equivalent conditions hold, one has moreover the following statements$:$
\begin{enumerate}[{\rm (i)}]
\item \label{Item1} for any compact set $K \subset \mathbb{R}^m,$ the set $f^{-1}(K - \mathbb{R}^m_+) \cap \Omega$ is compact$;$

\item \label{Item2} the sets $\mathrm{WSOL}$ and $\mathrm{WVAL}$ are closed$;$

\item \label{Item3} the following inclusion holds$:$
\begin{eqnarray*}
f(\Omega) &\subset& f({\rm ESOL}) + \mathbb{R}^m_+;
\end{eqnarray*}

\item \label{Item4} for any $\bar{x} \in \mathbb{R}^n,$ $\Omega(\bar{x}):= \mathrm{argmin}_{y \in \Omega} \phi({\bar{x}}, y)$ is a nonempty compact set and there exist constants $c > 0$ and $\alpha > 0$ such that
\begin{eqnarray*}
\phi({\bar{x}}, y) + \psi(\bar{x}) + \left[\phi({\bar{x}}, y) + \psi(\bar{x}) \right]^\alpha  &\ge& c \, \mathrm{dist}(y, \Omega(\bar{x}))
 \quad \textrm{ for all } \quad y \in \Omega;
\end{eqnarray*}

\item \label{Item5} for any $\bar{z} \in \mathrm{EVAL},$ there exist constants $c > 0$ and $\alpha > 0$ such that
\begin{eqnarray*}
\max_{i = 1, \ldots, m} \big(f_i(y) - \bar{z}_i \big) + \left[\max_{i = 1, \ldots, m} \big(f_i(y) - \bar {z}_i \big) \right]^\alpha 
&\ge& c\, \mathrm{dist}(y, \mathrm{ESOL}) \quad \textrm{ for all } \quad y \in \Omega;
\end{eqnarray*}

\item \label{Item6} for any $\bar{z} \in \mathrm{WVAL},$ there exist constants $c > 0$ and $\alpha > 0$ such that
\begin{eqnarray*}
\max_{i = 1, \ldots, m} \big(f_i(y) - \bar{z}_i \big) + \left[\max_{i = 1, \ldots, m} \big(f_i(y) - \bar {z}_i \big) \right]^\alpha 
&\ge& c\, \mathrm{dist}(y, \mathrm{WSOL}) \quad \textrm{ for all } \quad y \in \Omega;
\end{eqnarray*}

\item \label{Item7} for any $\bar{z} \in \mathbb{R}^m,$ there exist constants $c > 0, \alpha > 0,$ and $r > 0$ such that
\begin{eqnarray*}
{L}(z) &\subset& {L}(\bar{z}) + c\|z - \bar{z}\|^\alpha \mathbb{B} \quad \textrm{ for all } \quad z \in \mathbb{B}_r(\bar{z}),
\end{eqnarray*}
and hence in particular the set-valued mapping $L \colon \mathbb{R}^{m} \rightrightarrows \mathbb{R}^n$ is upper semicontinuous\footnote{A set-valued mapping ${L} \colon \mathbb{R}^{m} \rightrightarrows \mathbb{R}^n$ is said to be {\em upper semicontinuous} at $\bar{z} \in \mathbb{R}^{m}$ if for each open set $V \subset \mathbb{R}^{n}$ satisfying ${L} (\bar{z}) \subset V,$ there exists $r > 0$ such that ${L}(z) \subset V$ whenever $z \in \mathbb{B}_r(\bar{z}).$};

\item \label{Item8}  if $x_k \in \Omega$ is a sequence such that $f(x_k)$ converges to some value $\bar{z} \in \mathrm{EVAL}$ $($resp., $\bar{z} \in \mathrm{WVAL}),$ then the sequence $x_k$ is bounded, and all its cluster points lie in $\mathrm{ESOL}$ $($resp., $\mathrm{WSOL}).$
\end{enumerate}
\end{theorem}

\begin{proof}
We first show the equivalence between the conditions (c1)--(c7).

(c1) $\Rightarrow$ (c2) This is clear.

(c2) $\Rightarrow$ (c3) Assume that the set of all weakly efficient solutions to the problem~\eqref{MP} contains a point, say $\bar{x}.$ Then for $\bar{z} := f(\bar{x}) - (1, \ldots, 1) \in \mathbb{R}^m$ we have
\begin{eqnarray*}
f(\Omega) \cap \left(\bar z - \mathbb{R}_+^m\right)& = & \emptyset,
\end{eqnarray*}
which is the desired condition.

(c3) $\Rightarrow$ (c4) Take any $\bar{x} \in \mathbb{R}^n$ and consider the function $h \colon \mathbb{R}^n \to \mathbb{R}, x \mapsto \max_{i = 1, \ldots, m} f_i(x).$ Then $h$ is a continuous selection of the polynomials $f_1, \ldots, f_m$ (by definition), semialgebraic (by the Tarski--Seidenberg theorem) and bounded from below (by Theorem~\ref{Thm51}). In view of Theorem~\ref{Thm41}, there exist constants $c_1 > 0$ and $R > 0$ such that for all $y \in \Omega \setminus \mathbb{B}_R$ we have
\begin{eqnarray*}
h(y) &\ge& c_1 \|y\|.
\end{eqnarray*}
Observe that
\begin{eqnarray*}
\phi({\bar{x}}, y) &\ge& h(y) - \max_{i = 1, \ldots, m} f_i(\bar{x}).
\end{eqnarray*}
Therefore, for all $y \in \Omega$ we have
\begin{eqnarray*}
\phi({\bar{x}}, y)&\ge& c_1\|y\| + c_2,
\end{eqnarray*}
with $c_2$ being defined by
\begin{eqnarray*}
c_2 &:=& \min\{- \max_{i = 1, \ldots, m} f_i(\bar{x}), \min_{y \in \mathbb{B}_R} \phi({\bar{x}}, y) - c_1R\}.
\end{eqnarray*}

(c4) $\Rightarrow$ (c5) This is clear.

(c5) $\Rightarrow$ (c6)
Let $\bar{x} \in \mathbb{R}^n$ be such that the function $\phi({\bar{x}}, \cdot)$ is coercive on $\Omega,$ and take any $\bar{z} \in \mathbb{R}^m.$ We have
\begin{eqnarray*}
\{x \in \Omega \mid f(x) \le \bar{z}\} &\subset &
\{x \in \Omega \mid \phi({\bar{x}}, x) \le \max_{i = 1, \ldots, m} ( \bar{z}_i - f_i(\bar{x})) \}.
\end{eqnarray*}
The right side of the inclusion is a bounded set, hence so is the left side. This, in combination with the continuity of $f,$ gives the required conclusion.

(c6) $\Rightarrow$ (c7) Take any $\bar{x} \in \Omega.$ By assumption, the set
\begin{eqnarray*}
\{x \in \Omega \mid f(x) \le f(\bar{x}) \}
\end{eqnarray*}
is nonempty compact. Since $f$ is continuous, it follows that the set
\begin{eqnarray*}
\{f(x) \mid x \in \Omega \textrm{ and } f(x) \le f(\bar{x}) \}
\end{eqnarray*}
is nonempty compact, as required.

(c7) $\Rightarrow$ (c1) Assume that for some $\bar{z} \in \mathbb{R}^m,$ the set $f(\Omega) \cap \left(\bar z - \mathbb{R}_+^m\right)$ is nonempty and bounded. It turns out that the set ${L}(\bar{z})$ is nonempty and that the restriction of $\max_{i = 1, \ldots, m} f_i(x)$ on $\Omega$ is bounded from below, and so is coercive by Theorem~\ref{Thm41}.

Fix any $\lambda \in \mathrm{int} \mathbb{R}^m_+.$ We show that the problem
\begin{eqnarray} \label{Eqn7}
\inf_{x \in {L}(\bar{z})} \, \langle \lambda, f(x) \rangle &=& \inf_{z \in f(\Omega) \cap \left(\bar z - \mathbb{R}_+^m\right)} \, \langle \lambda, z \rangle \ > \ -\infty
\end{eqnarray}
has a global minimizer. Indeed, if this were not true, there would exist a sequence $x_k \in {L}(\bar{z})$ such that $\|x_k\| \to \infty$ and $\langle \lambda, f(x_k) \rangle \to \inf_{x \in {L}(\bar{z})}\,\langle \lambda, f(x) \rangle.$ But then the sequence $f(x_k)$ must be bounded, which contradicts the fact that the restriction of $\max_{i = 1, \ldots, m} f_i(x)$ on $\Omega$ is coercive. Therefore, the problem~\eqref{Eqn7} has a global minimizer, say $\bar{x}.$ We have that $\bar{x}$ is an efficient solution for the problem~\eqref{MP}; otherwise, there would exist $\bar{x}' \in \Omega$ with $f(\bar{x}') - f(\bar{x})\in -\mathbb{R}^m_+ \setminus \{0\}.$ Certainly, $f(\bar{x}') \le f(\bar{x})\le \bar{z},$ and so $\bar{x}' \in {L}(\bar{z});$ moreover, we have
\begin{eqnarray*}
\langle \lambda, f(\bar{x}') \rangle &<& \langle \lambda, f(\bar{x}) \rangle \ = \ \inf_{x \in {L}(\bar{z})} \, \langle \lambda, f(x) \rangle,
\end{eqnarray*}
which is impossible.

\medskip
Next we show the statements \eqref{Item1}--\eqref{Item8} under the equivalent conditions (c1)--(c7).

\medskip
\eqref{Item1} Let $K$ be a compact subset of $\mathbb{R}^m.$ Then there exists a point $\bar{z} \in \mathbb{R}^m$ such that $K \subset \bar{z} - \mathbb{R}^m_+.$ Clearly,
\begin{eqnarray*}
f^{-1}(K - \mathbb{R}^m_+) \cap \Omega &\subset& \{x \in \Omega \mid f(x) \le \bar{z}\},
\end{eqnarray*}
which, together with Condition~(c6), the continuity of $f$ and the closedness of $\Omega,$ gives the desired statement.

\eqref{Item2} Note that from Condition~(c5) and Lemma~\ref{Lem51} that the function $\psi$ is finite, locally Lipschitz and semialgebraic.
To see the closedness of $\mathrm{WSOL},$ let $x_k \in \mathrm{WSOL}$ be a sequence tending to $\bar{x}.$ Then $x_k \in \Omega$ and $\psi(x_k) = 0.$ Hence $\bar{x} \in \Omega$ (since the set $\Omega$ is closed) and $\psi(\bar{x}) = 0$ (since the function $\psi$ is continuous).
Therefore, $\bar{x} \in \mathrm{WSOL},$ as required.

To see the closedness of $\mathrm{WVAL},$ let $z_k \in \mathrm{WVAL}$ be a sequence tending to $\bar{z};$ we must show that $\bar{z} \in \mathrm{WVAL}.$ Indeed, for each $k,$ there exists $x_k \in f^{-1}({z_k}) \cap \Omega.$ Take a compact neighborhood $K \subset \mathbb{R}^m$ of $\bar{z}.$ Then for all $k$ sufficiently large, $x_k$ belongs to $f^{-1}(K - \mathbb{R}^m_+) \cap \Omega,$ a compact set (by~\eqref{Item1}). Hence  the sequence $x_k$ has a cluster point, say $\bar{x}.$ Certainly, $\bar{x} \in \Omega.$ Note that $f(x_k) = z_k \to \bar{z}$ and $\psi(x_k) = 0.$ By continuity, $f(\bar{x}) = \bar{z}$ and $\psi(\bar{x}) = 0.$ Therefore, $\bar{x}\in \mathrm{WSOL},$ and so $\bar{z} \in \mathrm{WVAL},$ as required.

\eqref{Item3} Take any $\bar{z} \in f(\Omega).$ By Condition~(c6),
$\{x \in \Omega \mid f(x) \le \bar{z}\}$ is a nonempty compact set. From the proof of the implication (c7)~$\Rightarrow$~(c1), we can find a point $\bar{x} \in \mathrm{ESOL}$ with $f(\bar{x}) \le \bar{z}.$ Then $\bar{z} \in f(\mathrm{ESOL}) + \mathbb{R}^m_+.$ This is true for arbitrary $\bar{z} \in f(\Omega),$ so we conclude that
\begin{eqnarray*}
f(\Omega) &\subset& f({\rm ESOL}) + \mathbb{R}^m_+ ,
\end{eqnarray*}
as required.

\eqref{Item4} Let $\bar{x}$ be any point in $\mathbb{R}^n.$ Condition~(c4) ensures the existence of constants $c_1 > 0$ and $c_2 \in \mathbb{R}$ such that for all $y \in \Omega,$
\begin{eqnarray*}
\phi({\bar{x}}, y) &\ge& c_1\|y\| + c_2.
\end{eqnarray*}
It turns out that $\psi(\bar{x}) < \infty$ and that $\Omega(\bar{x}) $ is a nonempty compact set, and so $M := \sup_{x \in \Omega(\bar{x})} \|x\| < \infty.$ Moreover, it is easy to check that for all $y \in \mathbb{R}^n,$
\begin{eqnarray*}
\|y\| &\ge& \mathrm{dist}(y, \Omega(\bar{x})) - M.
\end{eqnarray*}
Choose $R$ large enough so that
\begin{eqnarray*}
- \frac{c_1}{2} M + \frac{c_1}{2}R  + c_2 + \psi(\bar{x}) &\ge& 0.
\end{eqnarray*}
Then we have for any $y \in \Omega \setminus \mathbb{B}_R$ that
\begin{eqnarray*}
\phi({\bar{x}}, y) + \psi(\bar{x}) 
&\ge& c_1\|y\| + c_2 + \psi(\bar{x})  \\
&\ge& \frac{c_1}{2} \big (\mathrm{dist}(y, \Omega(\bar{x}) ) - M \big) + \frac{c_1}{2} \|y\| + c_2 + \psi(\bar{x})  \\
&\ge& \frac{c_1}{2} \, \mathrm{dist}(y, \Omega(\bar{x})) - \frac{c_1}{2} M + \frac{c_1}{2} R + c_2 + \psi(\bar{x})  \\
&\ge& \frac{c_1}{2} \, \mathrm{dist}(y, \Omega(\bar{x}) ).
\end{eqnarray*}
On the other hand, by the classical \L ojasiewicz inequality (see, for example, \cite[Corollary~1.2]{HaHV2017}), there exist constants $c_3 > 0$ and $\alpha > 0$ such that for all $y \in \Omega \cap \mathbb{B}_R,$ 
\begin{eqnarray*}
\left[\phi({\bar{x}}, y) + \psi(\bar{x}) \right]^\alpha &\ge& c_3\, \mathrm{dist}(y, \Omega(\bar{x})).
\end{eqnarray*}
Letting $c := \min\{\frac{c_1}{2}, c_3\} > 0,$ we get for all $y \in \Omega$ that
\begin{eqnarray*}
\phi({\bar{x}}, y) + \psi(\bar{x}) + \left[\phi({\bar{x}}, y) + \psi(\bar{x}) \right]^\alpha  &\ge& c \, \mathrm{dist}(y, \Omega(\bar{x})).
\end{eqnarray*}

\eqref{Item5} Take any $\bar{z} \in \mathrm{EVAL}$ and fix $\bar{x} \in f^{-1}(\bar{z}) \cap \Omega \subset \mathrm{ESOL}.$ Then $\psi(\bar{x}) = 0$ and $\Omega(\bar{x}) \subset f^{-1}(\bar{z}) \cap \Omega $ (by definition). In combination with \eqref{Item4}, the desired statement follows.

\eqref{Item6} Take any $\bar{z} \in \mathrm{WVAL}$ and fix $\bar{x} \in f^{-1}(\bar{z}) \cap \Omega.$ Then $\psi(\bar{x}) = 0$ and $\Omega(\bar{x}) \subset \mathrm{WSOL}$ (by Lemma~\ref{Lem51}). In combination with \eqref{Item4}, the desired statement follows.

\eqref{Item7} Fix any $\bar{z} \in \mathbb{R}^m.$ If $L(\bar{z}) = \emptyset,$ then there exists $r > 0$ such that $L(z) = \emptyset$ for all $z \in \mathbb{B}_r(\bar{z}).$ (If not, we would have sequences $z_k \to \bar{z}$ and $x_k \in L(z_k).$ Certainly the sequence $x_{k}$ is bounded, and hence has a cluster point, which clearly lies in ${L}(\bar{z}),$ a contradiction). The desired statement holds trivially in this case.

Assume that $L(\bar{z}) \ne \emptyset$ and take any $r > 0.$ The closed ball $K := \mathbb{B}_r(\bar{z})$ is compact. By \eqref{Item1}, the set $f^{-1}(K - \mathbb{R}^m_+) \cap \Omega$ is compact, and so is every set ${L}(z) := f^{-1}(z - \mathbb{R}^m_+) \cap \Omega,$ $z \in K.$  Define the function $g \colon K \to \mathbb{R}$ by
\begin{eqnarray*}
g(z) &:=& 
\begin{cases}
\sup_{x \in {L}(z)} \mathrm{dist}(x, {L}(\bar{z})) & \textrm{ if } L(z) \ne \emptyset, \\
0 & \textrm{ otherwise,}
\end{cases}
\end{eqnarray*}
which is well-defined, semialgebraic (by the Tarski--Seidenberg theorem), and satisfies
\begin{eqnarray*}
0 &\le& g(z) \ \le \ M \ := \ \sup_{x \in f^{-1}(K - \mathbb{R}^m_+) \cap \Omega} \mathrm{dist}(x, {L}(\bar{z})) \ < \ \infty.
\end{eqnarray*}
Furthermore, for every sequence $z_k \in K$ converging to $\bar{z},$ the corresponding sequence $g(z_k)$ converges to $0.$ (Indeed, 
let $z_{k_l}$ be a subsequence of $z_k$ such that either $L(z_{k_l}) = \emptyset$  for all $l$ or, quite the opposite, $L(z_{k_l}) \ne \emptyset$ 
for all $l.$ In the first case, $g(z_{k_l}) = 0$ for all $l.$ In the second case, there exists a sequence $x_{k_l} \in {L}(z_{k_l}) \subset f^{-1}(K - \mathbb{R}^m_+) \cap \Omega$ such that $g(z_{k_l}) = \mathrm{dist}(x_{k_l}, {L}(\bar{z}));$ certainly the sequence $x_{k_l}$ is bounded, and all its cluster points lie in ${L}(\bar{z}),$ which imply that $g(z_{k_l}) \to 0$ as $l \to \infty.$) 

We next consider the function
$$\mu \colon [0, r] \to \mathbb{R}, \quad t\mapsto \sup_{\|z - \bar{z} \| = t}g(z),$$
which is well-defined, semialgebraic (by the Tarski--Seidenberg Theorem), and satisfies $\mu(t) \to 0$ as $t \to 0.$
By Lemmas~\ref{MonotonicityLemma} and \ref{GrowthDichotomyLemma}, there exist positive constants $c', \alpha$ and $\epsilon$ such that for all $t \in [0, \epsilon],$
\begin{eqnarray*}
\mu(t) &\le& c' t^\alpha.
\end{eqnarray*}
It turns out that if $\|z - \bar{z}\| \le \epsilon,$ then
\begin{eqnarray*}
g(z) &\le& \mu(\|z - \bar{z}\|) \ \le \ c'\|z - \bar{z}\|^\alpha.
\end{eqnarray*}
On the other hand, if $\epsilon \le \|z - \bar{z}\| \le r,$ then
\begin{eqnarray*}
g(z) &\le& M \ = \ \frac{M}{\epsilon^\alpha} \epsilon^\alpha \ \le \ \frac{M}{\epsilon^\alpha}  \|z - \bar{z}\|^\alpha.
\end{eqnarray*}
Therefore, letting $c := \max\{c', \frac{M}{\epsilon^\alpha} \},$ we obtain for all $z \in K,$
\begin{eqnarray*}
g(z) &\le& c\|z - \bar{z}\|^\alpha,
\end{eqnarray*}
which yields the desired conclusion
\begin{eqnarray*}
{L}(z) &\subset& {L}(\bar{z}) + c\|z - \bar{z}\|^\alpha \mathbb{B}.
\end{eqnarray*}

\eqref{Item8} Let $x_k$ be a sequence in $\Omega$ such that $f(x_k)$ tends to some value $\bar{z} \in \mathbb{R}^m.$
By \eqref{Item1}, the sequence $x_k$ is bounded, hence by \eqref{Item7}, all its cluster points lie in ${L}(\bar{z}).$
The proof is complete by noting that if $\bar{z} \in \mathrm{EVAL}$ (resp., $\bar{z} \in \mathrm{WVAL}$) then
${L}(\bar{z}) \subset \mathrm{ESOL}$ (resp., ${L}(\bar{z}) \subset \mathrm{WSOL}$).
\end{proof}

\begin{remark}{\rm
(i) The condition~(c6) implies that the set $f(\Omega)$ is {\em $\mathbb{R}^m_+$-compact} in the sense that for any $\bar{z} \in f(\Omega),$ the set  $f(\Omega) \cap (\bar{z} - \mathbb{R}^m_+)$ is compact; for related results about this notion, see \cite{Ehrgott2005, Luc1989}.

(ii) The inclusion in \eqref{Item3} is called the {\em domination property} (or sometimes called {\em external stability}) for the set $f(\Omega)$ (see
\cite[Definition~4.1]{Luc1989} and \cite[Definition~2.20]{Ehrgott2005}).

(iii) In the scalar optimization case (i.e., $m = 1$), $\mathrm{ESOL} = \mathrm{WSOL},$ $\mathrm{EVAL} = \mathrm{WVAL} = \{\inf_{x \in \Omega} f(x)\},$ and the inequalities in \eqref{Item5} and \eqref{Item6} read as follows:
\begin{eqnarray*}
f(y) - \inf_{x \in \Omega} f(x) + \left [ f(y) - \inf_{x \in \Omega} f(x) \right]^\alpha &\ge& c\, \mathrm{dist}(y, \mathrm{ESOL})
\quad \textrm{ for all } \quad y \in \Omega.
\end{eqnarray*}

(iv) When the constraint set $\Omega$ is defined by polynomial equations and inequalities, by \cite[Theorem~3.3]{HaHV2017}), the exponent $\alpha$ can be explicitly estimated from the data of the problem. Since we do not use this fact, we leave the details for the interested reader.
}\end{remark}

The rest of this section is devoted to the existence of properly efficient solutions. To see this, recall that the {\em recession cone} $[Z]^{\oplus}$ of a set $Z \subset \mathbb{R}^m$ is defined by
\begin{eqnarray*}
[Z]^{\oplus} &:=& \big\{d\in\mathbb{R}^m\,\big|\,\exists\,\{(t_k, z_k)\}\subset\mathbb{R}_+\times Z \;\text{ such that }\;t_k \to 0 \;\textrm{ and }\; t_k z_k\to d\big\}.
\end{eqnarray*}
Theorem~\ref{Thm52} allows us to prove the following result, which provides necessary and sufficient conditions for the existence of properly efficient solutions.

\begin{corollary} \label{Corollary51}
Assume that the mapping $f$ is convenient and the restriction $f|_\Omega$ is nondegenerate at infinity. Then 
the problem~\eqref{MP} admits a properly efficient solution if and only if
\begin{eqnarray}\label{Eqn10}
[f(\Omega)]^{\oplus}\cap(-\mathbb{R}^m_+) &=&\{0\}.
\end{eqnarray}
\end{corollary}

\begin{proof}
Assume that the problem~\eqref{MP} admits a properly efficient solution. In view of \cite[Theorem~5.2]{Henig1982} (see also \cite[Proposition~2.1]{Luc1990}), we have
\begin{eqnarray*}
[f(\Omega)+\mathbb{R}^m_+]^{\oplus}\cap(-\mathbb{R}^m_+)&=&\{0\},
\end{eqnarray*}
from which the equation~\eqref{Eqn10} is evident.

Conversely, assume that \eqref{Eqn10} holds. Then the restriction $f|_\Omega$ has an empty section. (Indeed, if this were not true, there
would exist a sequence $x_k \in \Omega$ such that $\max_i f_i(x_k) \to -\infty$ as $k \to \infty.$ The sequence $f(x_k)/\|
f(x_k)\|$ has at least one cluster point, say $d.$ Certainly, $d \ne 0$ and $d \in [f(\Omega)]^{\oplus}\cap(-\mathbb{R}^m_+) ,$ a contradiction.)
Therefore, the equivalent conditions~(c1)--(c7) in Theorem~\ref{Thm52} hold.

We next show that $f(\Omega) + \mathbb{R}^m_+$ is a closed set. To see this, pick an arbitrary sequence $z_k\in f(\Omega)+\mathbb{R}^m_+$ that converges to some $z \in \mathbb{R}^m$ and find sequences $x_k\in \Omega$ and $d_k \in \mathbb{R}^m_+$ such that $z_k = f(x_k) + d_k$ for all $k.$ Since the sequence $z_k$ is convergent, there is $\bar{z} \in \mathbb{R}^m$ with $z_k \le \bar{z}$ for all $k.$ Certainly, $f(x_k) \le \bar{z},$ and hence $x_k \in  L({\bar{z}}),$ which is a compact set by Condition~(c6) in Theorem~\ref{Thm52}.
Consequently, there is a subsequence of $x_k$ converging to some $x \in L({\bar{z}}).$ It turns out that the corresponding subsequence of $d_k$ is also convergent to some $d\in\mathbb{R}^m_+.$ Therefore,
\begin{eqnarray*}
z &=& f(x) + d \ \in \ f(L({\bar{z}})) + \mathbb{R}^m_+ \ \subset \ f(\Omega)+\mathbb{R}^m_+,
\end{eqnarray*}
which proves the closedness of $f(\Omega) + \mathbb{R}^m_+.$ This, in combining with \eqref{Eqn10}, \cite[Theorem~3.2]{Benson1979} and \cite[Theorems~2.1~and~5.1]{Henig1982}, then yields the desired conclusion.
\end{proof}

\begin{remark}{\rm
It is well known that, in the setting under consideration, the concept of properly efficient solutions in Geoffrion's sense agrees with those of properly efficient solutions in the senses of Benson~\cite{Benson1979} and Henig~\cite{Henig1982}, and is stronger than that of properly efficient solutions in the senses of Borwein~\cite{Borwein1977}; see \cite{Ehrgott2005, Luc1989}.
}\end{remark}

\section{Examples and remarks} \label{Section6}

In this section we provide some examples and remarks concerning the results presented in the previous section.
Let us start with the following example which is an illustration of Theorem~\ref{Thm52}.

\begin{example}\label{Example61}{\rm
Let $\Omega := \mathbb{R}^2$ and consider the polynomial mapping
\begin{eqnarray*}
f := (f_1, f_2) \colon \mathbb{R}^2 \to \mathbb{R}^2, && (x_1, x_2) \mapsto \left(3x_1 + 4 x_2,\ 25(x_1^2 + x_2^2)\right).
\end{eqnarray*}
By definition, $\Gamma(f_1) = {\rm conv}\left\{(1,0),(0,1)\right\}$ and $\Gamma(f_2) = {\rm conv}\left\{(2,0),(0,2)\right\},$ and hence $f$ is convenient (see {\sc Figure~\ref{fig0512-1}}).

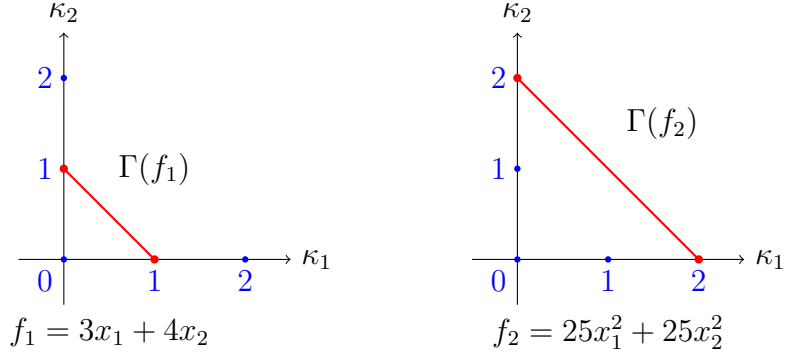
\begin{figure}[h!]
	\centering
\begin{tikzpicture}[scale=1.2]

\begin{scope}[shift={(-2.5,0)}]
  \draw[->] (-0.5,0) -- (2.5,0) node[right] {$\kappa_1$};
  \draw[->] (0,-0.5) -- (0,2.5) node[above] {$\kappa_2$};
  \fill[blue] (2,0) circle (1.0pt) node[below] {$2$};
  \fill[blue] (1,0) circle (1.0pt) node[below] {$1$};
  \fill[blue] (0,1) circle (1.0pt) node[left] {$1$};
  \fill[blue] (0,2) circle (1.0pt) node[left] {$2$};
  \fill[blue] (0,0) circle (1.0pt) node[below left] {$0$};
  \draw[thick, red] (1,0) -- (0,1);
  \fill[red] (1,0) circle (1.3pt) node[right] {$\quad$};
  \fill[red] (0,1) circle (1.3pt) node[right] {$\quad$};
  \node at (1,1) {$\Gamma(f_1)$};
  \node at (0.5,-0.8) {$f_1 = 3x_1 + 4x_2$};
\end{scope}

\begin{scope}[shift={(2.5,0)}]
  \draw[->] (-0.5,0) -- (2.5,0) node[right] {$\kappa_1$};
  \draw[->] (0,-0.5) -- (0,2.5) node[above] {$\kappa_2$};
  \fill[blue] (2,0) circle (1.0pt) node[below] {$2$};
  \fill[blue] (0,2) circle (1.0pt) node[left] {$2$};
  \fill[blue] (1,0) circle (1.0pt) node[below] {$1$};
  \fill[blue] (0,1) circle (1.0pt) node[left] {$1$};
  \fill[blue] (0,0) circle (1.0pt) node[below left] {$0$};
  \draw[thick, red] (2,0) -- (0,2);
  \fill[red] (2,0) circle (1.3pt) node[right] {$\quad$};
  \fill[red] (0,2) circle (1.3pt) node[right] {$\quad$};
  \node at (1.6,1.5) {$\Gamma(f_2)$};
  \node at (1,-0.8) {$f_2 = 25x_1^2 + 25x_2^2$};
\end{scope}
\end{tikzpicture}
\caption{\small The Newton polyhedra of $f_1$ and $f_2$ in Example~\ref{Example61}.}
\label{fig0512-1}
\end{figure}

We next verify that the mapping $f$ is nondegenerate at infinity. For any $q = (q_1, q_2) \in \mathbb{R}^2$ with $\min\{q_1, q_2\} < 0,$ the following three cases happen.
\begin{itemize}
\item If $q_1 > q_2,$ then
\begin{itemize}
\item[$\circ$] $d_1 = d(q, \Gamma(f_1)) = q_2 < 0$ and $\Delta_1 = \Delta(q, \Gamma(f_1)) = \{(0,1)\} \subset \mathbb{R}^{\{2\}};$
\item[$\circ$] $d_2 = d(q, \Gamma(f_2)) = 2q_2 < 0$ and $\Delta_2 = \Delta(q, \Gamma(f_2)) = \{(0,2)\}\subset \mathbb{R}^{\{2\}};$
\item[$\circ$] $f_{1, \Delta_1}(x_1, x_2) = 4 x_2$ and $f_{2, \Delta_2}(x_1, x_2) = 25 x_2^2.$
\end{itemize}
\item If $q_1 < q_2,$ then
\begin{itemize}
\item[$\circ$] $d_1 = d(q, \Gamma(f_1)) = q_1 < 0$ and $\Delta_1 = \Delta(q, \Gamma(f_1)) = \{(1,0)\}\subset \mathbb{R}^{\{1\}};$
\item[$\circ$] $d_2 = d(q, \Gamma(f_2)) = 2q_1 < 0$ and $\Delta_2 = \Delta(q, \Gamma(f_2)) = \{(2,0)\}\subset \mathbb{R}^{\{1\}};$
\item[$\circ$] $f_{1, \Delta_1}(x_1, x_2) = 3 x_1$ and $f_{2, \Delta_2}(x_1, x_2) = 25 x_1^2.$
\end{itemize}
\item If $q_1 = q_2,$ then
\begin{itemize}
\item[$\circ$] $d_1 = d(q, \Gamma(f_1)) = q_1 = q_2 < 0$ and $\Delta_1 = \Delta(q, \Gamma(f_1)) = {\rm conv}\{(1,0), (0,1)\}\subset \mathbb{R}^{\{1,2\}};$
\item[$\circ$] $d_2 = d(q, \Gamma(f_2)) = 2q_1 = 2 q_2 < 0$ and $\Delta_2 = \Delta(q, \Gamma(f_2)) = {\rm conv}\{(2,0), (0,2)\}\subset \mathbb{R}^{\{1,2\}};$
\item[$\circ$] $f_{1, \Delta_1}(x_1, x_2) = 3 x_1 + 4 x_2$ and $f_{2, \Delta_2}(x_1, x_2) = 25 x_1^2 + 25 x_2^2.$
\end{itemize}
\end{itemize}
In each case, we can see that for any nonempty sets $I \subset \{1, 2\}$ and $J \subset \{1, 2\},$ there are no $(c_1, c_2) \in (\mathbb{R}^*_+)^I$ and $(x_1^*, x_2^*) \in (\mathbb{R}^*)^J$ satisfying the conditions~(c1)-(c3) in Definition~\ref{nondegeneracy}. Therefore, $f$ is nondegenerate at infinity.
On the other hand, it is not hard to see that
\begin{eqnarray*}
f(\mathbb{R}^2) &=& \left\{(z_1, z_2) \in \mathbb{R}^2 \mid z_2 \geq z_1^2\right\}
\end{eqnarray*}
(see {\sc Figure~\ref{fig0512-2}}(b)); in particular, the mapping $f$ has nonempty bounded sections as well as empty sections.
By Theorem~\ref{Thm52}, the sets $\mathrm{ESOL}$ and $\mathrm{WSOL}$ are nonempty.  Indeed, it is easy to see that
\begin{eqnarray*}
{\rm ESOL} \ = \ {\rm WSOL} &=& \left\{(x_1, x_2) \mid x_2 = \tfrac{4}{3}x_1, x_1 \leq 0\right\}, \\
{\rm EVAL} \ = \ {\rm WVAL} &=& \left\{(z_1, z_2) \mid z_2 = z_1^2, z_1 \leq 0 \right\},
\end{eqnarray*}
(see {\sc Figure~\ref{fig0512-2}}). Also, the achievement function is given by
\begin{eqnarray*}
\psi(x_1, x_2) =
\begin{cases}
25(x_1^2+x_2^2) & \text{if } \sqrt{\left(x_1-\frac{3}{50}\right)^2+\left(x_2-\frac{4}{50}\right)^2}\le\frac{1}{10},\\[10pt]
3x_1+4x_2 - \frac{1}{2} + 5 \sqrt{\left(x_1-\frac{3}{50}\right)^2 + \left(x_2-\frac{4}{50}\right)^2} & \text{otherwise}.
\end{cases}
\end{eqnarray*}
By definition, we have for any $z := (z_1, z_2) \in \mathbb{R}^2$ that
\begin{eqnarray*}
L(z) &=& f^{-1}\left(z - \mathbb{R}^2_+\right) \cap \mathbb{R}^2 \\
&=& \begin{cases}
\emptyset & \textrm{ if } {z}_2 < 0, \\
\left\{(x_1, x_2) \in \mathbb{R}^2 \mid 3x_1 + 4x_2 \leq z_1,\ 25(x_1^2 + x_2^2) \leq z_2 \right\} & \textrm{ otherwise.}
\end{cases}
\end{eqnarray*}
A simple (but tedious) calculation shows that for any $\bar{z} \in \mathrm{WVAL} = \mathrm{EVAL},$ there exists a constant $c > 0$ such that
\begin{eqnarray*}
\max_{i = 1, 2} \big(f_i(y) - \bar{z}_i \big) + \left[\max_{i = 1, 2} \big(f_i(y) - \bar{z}_i \big) \right]^{\frac{1}{2}} &\ge& c\, \mathrm{dist}(y, \mathrm{WSOL}) \ = \ c\, \mathrm{dist}(y, \mathrm{ESOL})
\end{eqnarray*}
for all $y \in \Omega.$
Finally, for any $\bar{z} := (\bar{z}_1, \bar{z}_2) \in \mathbb{R}^2,$ there exist constants $c > 0$ and $r > 0$ such that
\begin{eqnarray*}
L(z) &\subset& L(\bar{z}) + c\|z - \bar{z}\|^\alpha \mathbb{B} \quad \textrm{ for all } \quad z \in \mathbb{B}_r(\bar{z}),
\end{eqnarray*}
where
\begin{eqnarray*}
\alpha :=
\begin{cases}
\frac{1}{2} & \textrm{ if } \bar{z}_2 = \bar{z}_1^2, \\
1 & \textrm{ otherwise.}
\end{cases}
\end{eqnarray*}

\begin{figure}[h!]
	\centering
\begin{tikzpicture}[scale=0.9]
\begin{scope}[shift={(-4,0)}]
  \draw[->] (-3.1,0) -- (3.1,0) node[right] {$x_1$};
  \draw[->] (0,-1.5) -- (0,4.2) node[above] {$x_2$};
  \draw[thick, gray] plot[domain=0:2, samples=100] (\x, {4/3*\x})
    node[right] {};
  \draw[very thick, red] plot[domain=-1.1:0, samples=100] (\x, {4/3*\x})
    node[right] {};
  \node at (-1.5,-0.6) {ESOL};
  \node at (2,1.1) {$x_2 = \frac{4}{3}x_1$};
  \node at (0,-2) {${\rm (a)}$};
\end{scope}

\begin{scope}[shift={(4,0)}]
  \draw[->] (-3.1,0) -- (3.1,0) node[right] {$z_1$};
  \draw[->] (0,-1.5) -- (0,4.2) node[above] {$z_2$};
  \fill[gray!30, opacity=0.5]
    plot[domain=-2:2, samples=100] (\x, {\x*\x}) --
    (2,4.0) -- (-2,4.0) -- cycle;
  \draw[thick, gray] plot[domain=0:2, samples=100] (\x, {\x*\x})
    node[right] {};
  \draw[very thick, blue] plot[domain=-2:0, samples=100] (\x, {\x*\x})
    node[right] {};
  \node at (-0.8,2.8) {$f(\mathbb{R}^2)$};
  \node at (-1.9,0.9) {EVAL};
  \node at (2.2,1.6) {$z_2 = z_1^2$};
  \node at (0,-2) {${\rm (b)}$};
\end{scope}
\end{tikzpicture}
\caption{\small (a) The set of efficient solutions (red color); \ (b) the image $f(\mathbb{R}^2)$ (gray shadow) and the set of efficient values (blue color) in Example~\ref{Example61}.}
\label{fig0512-2}
\end{figure}
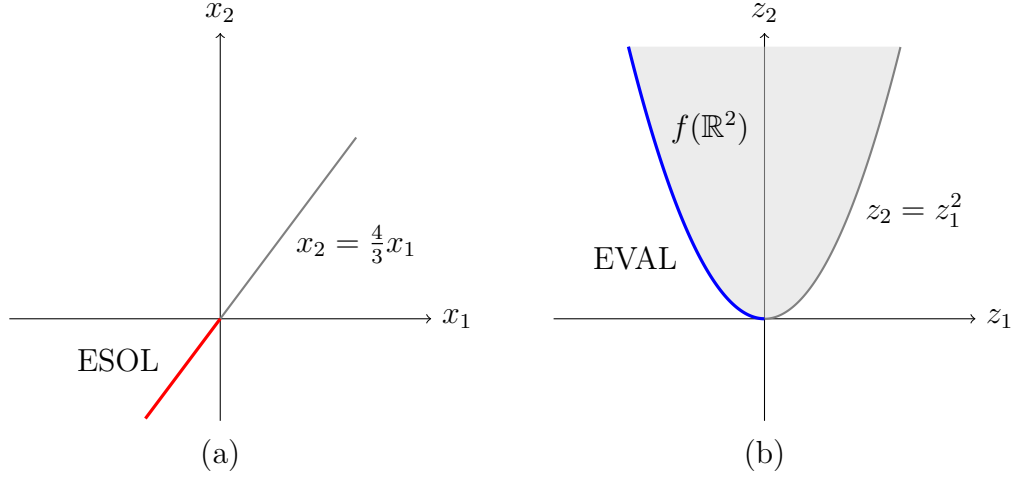
}\end{example}

The following example shows that the sets $\mathrm{ESOL}$ and $\mathrm{EVAL}$ may not be closed even when the assumptions of Theorem~\ref{Thm52} are satisfied. It also shows that, in general, the uniform bounded weak sharp minima does not hold.

\begin{example}\label{Example62}{\rm
Let $\Omega := \{(x_1, x_2) \in \mathbb{R}^2 \mid x_1^3 + x_2^3 - x_1 + x_2 \geq 0 \} \cup \{(1, -1)\}$ and consider the polynomial mapping
\begin{eqnarray*}
f \colon \mathbb{R}^2 \to \mathbb{R}^2, && (x_1, x_2) \mapsto \left(x_1^3 + x_2, -x_1 + x_2^3\right).
\end{eqnarray*}
By definition, we have
$C(\infty; \Omega) = \{(-t ,t) \mid t \in \mathbb{R}\},$ $\Gamma(f_1) = {\rm conv}\{(3,0),(0,1)\}$ and $\Gamma(f_2) = {\rm conv}\{(1,0),(0,3)\}$
(see {\sc Figure~\ref{fig0609-1}}). Then it is not hard to check that the mapping $f$ is convenient and the restriction $f|_\Omega$ is nondegenerate at infinity. Moreover, a direct calculation shows that
\begin{eqnarray*}
f(\Omega) &=& \left\{(z_1, z_2) \in \mathbb{R}^2 \mid z_1 + z_2 \geq 0 \right\} \cup \{(0, -2)\}
\end{eqnarray*}
(see {\sc Figure~\ref{fig0609-2}}(b)); particularly, the restriction $f|_{\Omega}$ has nonempty bounded sections as well as empty sections.
\begin{figure}[h!]
	\centering
\begin{tikzpicture}[scale=1]
\begin{scope}[shift={(-3.5,0)}]
  \draw[->] (-0.5,0) -- (3.5,0) node[right] {$\kappa_1$};
  \draw[->] (0,-0.5) -- (0,3.5) node[above] {$\kappa_2$};
  \fill[blue] (3,0) circle (1.0pt) node[below] {$3$};
  \fill[blue] (2,0) circle (1.0pt) node[below] {$2$};
  \fill[blue] (1,0) circle (1.0pt) node[below] {$1$};
  \fill[blue] (0,1) circle (1.0pt) node[left] {$1$};
  \fill[blue] (0,2) circle (1.0pt) node[left] {$2$};
  \fill[blue] (0,3) circle (1.0pt) node[left] {$3$};
  \fill[blue] (0,0) circle (1.0pt) node[below left] {$0$};
  \draw[thick, red] (3,0) -- (0,1);
  \fill[red] (0,1) circle (1.3pt) node[right] {$\quad$};
  \fill[red] (3,0) circle (1.3pt) node[right] {$\quad$};
  \node at (2,1) {$\Gamma(f_1)$};
  \node at (1.5,-0.8) {${\rm (a)}$};
\end{scope}

\begin{scope}[shift={(3.5,0)}]
  \draw[->] (-0.5,0) -- (3.5,0) node[right] {$\kappa_1$};
  \draw[->] (0,-0.5) -- (0,3.5) node[above] {$\kappa_2$};
  \fill[blue] (3,0) circle (1.0pt) node[below] {$3$};
  \fill[blue] (2,0) circle (1.0pt) node[below] {$2$};
  \fill[blue] (1,0) circle (1.0pt) node[below] {$1$};
  \fill[blue] (0,1) circle (1.0pt) node[left] {$1$};
  \fill[blue] (0,2) circle (1.0pt) node[left] {$2$};
  \fill[blue] (0,3) circle (1.0pt) node[left] {$3$};
  \fill[blue] (0,0) circle (1.0pt) node[below left] {$0$};
  \draw[thick, red] (1,0) -- (0,3);
  \fill[red] (0,3) circle (1.3pt) node[right] {$\quad$};
  \fill[red] (1,0) circle (1.3pt) node[right] {$\quad$};
  \node at (1.5,1) {$\Gamma(f_2)$};
  \node at (1.5,-0.8) {${\rm (b)}$};
\end{scope}
\end{tikzpicture}
\caption{\small (a) The Newton polyhedron of $f_1;$ \ (b) the Newton polyhedron of $f_2$ in Example~\ref{Example62}}
\label{fig0609-1}
\end{figure}
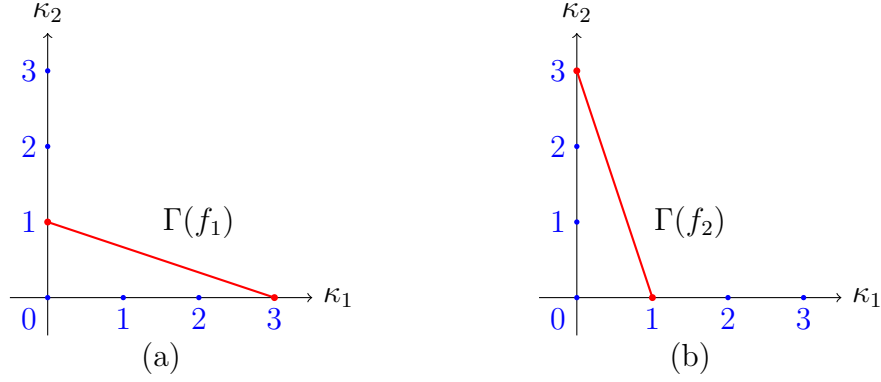

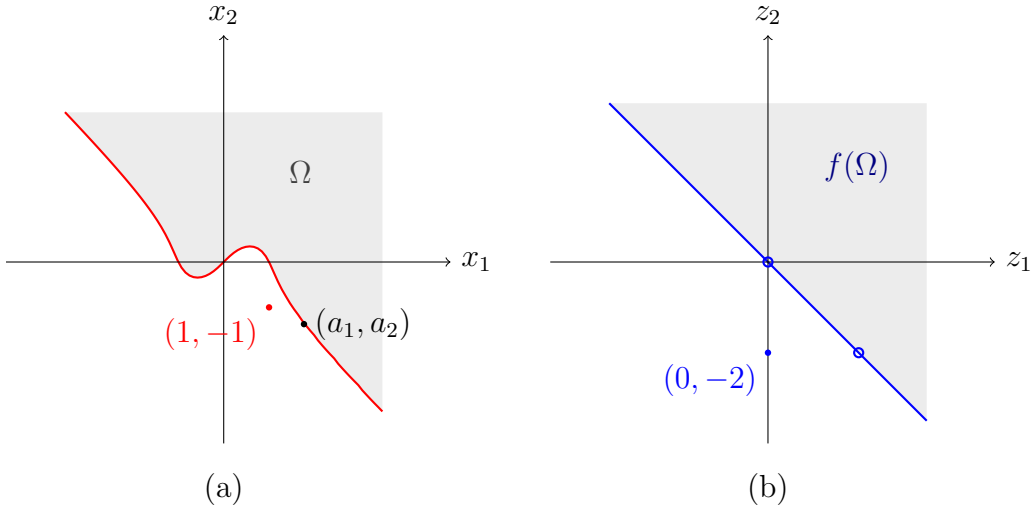
\begin{figure}[h!]
	\centering
\begin{tikzpicture}[scale=0.6]
\begin{scope}[shift={(-6,0)}]
  \fill[gray!15]
    plot[domain=-3.5:3.5, samples=80] (\x, {(2/sqrt(3)) * sinh((1/3) * ln(((3*sqrt(3)/2)*(\x - \x^3)) + sqrt( ((3*sqrt(3)/2)*(\x - \x^3))^2 + 1 )))}) -- (3.5,3.3) -- (-3.5,3.3) -- cycle;
  \draw[thick, red] plot[domain=-3.5:3.5, samples=80] (\x, {(2/sqrt(3)) * sinh((1/3) * ln(((3*sqrt(3)/2)*(\x - \x^3)) + sqrt( ((3*sqrt(3)/2)*(\x - \x^3))^2 + 1 )))})
    node[right] {};
  \draw[->] (-4.8,0) -- (5,0) node[right] {$x_1$};
  \draw[->] (0,-4) -- (0,5) node[above] {$x_2$};
  \node[above right, gray!50!black] at (1.2,1.5) {$\Omega$};
  \fill[red] (1,-1) circle (2pt) node[below left] {$(1,-1)$};
\fill[black] (1.77,-1.37) circle (2pt) node[right] {$(a_1,a_2)$};
  \node at (0,-5) {${\rm (a)}$};
\end{scope}

\begin{scope}[shift={(6,0)}]
  \fill[gray!15] (-3.5,3.5) -- (3.5,3.5) -- (3.5,-3.5) -- cycle;
  \draw[thick,blue] (-3.5,3.5) -- (3.5,-3.5) node[below right]{};
  \draw[->] (-4.8,0) -- (5,0) node[right] {$z_1$};
  \draw[->] (0,-4) -- (0,5) node[above] {$z_2$};
  \node[above right, blue!50!black] at (1,1.5) {$f(\Omega)$};
\draw[blue, thick] (0,0) circle (0.1);
\draw[blue, thick] (2,-2) circle (0.1);
  \fill[blue] (0,-2) circle (2pt) node[below left] {$(0,-2)$};
  \node at (0,-5) {${\rm (b)}$};
\end{scope}
\end{tikzpicture}
\caption{\small (a) The constraint set $\Omega$ including the point $(1, -1);$ \ (b) the image $f(\Omega)$ including the point $(0, -2)$ in Example~\ref{Example62}.}
\label{fig0609-2}
\end{figure}
By Theorem~\ref{Thm52}, the sets $\mathrm{ESOL}$ and $\mathrm{WSOL}$ are nonempty. Indeed, let $(a_1, a_2) \in \mathbb{R}^2$ be the unique solution of the system
\begin{eqnarray*}
x_1^3 + x_2 =  2, &\ &  -x_1 + x_2^3 =  -2.
\end{eqnarray*}
Then it is easy to check that
\begin{eqnarray*}
{\rm ESOL} &=& \left\{(x_1, x_2) \in \mathbb{R}^2 \mid x_1^3 + x_2^3 - x_1 + x_2 = 0,\ x_1 (x_1 - a_1) > 0\right\} \cup \left\{(1, -1)\right\}, \\
{\rm EVAL} &=& \left\{(z_1, z_2) \in \mathbb{R}^2 \mid z_1 + z_2 = 0, z_1 (z_1 - 2) > 0 \right\} \cup \left\{(0, -2)\right\},
\end{eqnarray*}
which are not closed, while the sets
\begin{eqnarray*}
{\rm WSOL} &=& \left\{(x_1, x_2) \in \mathbb{R}^2 \mid x_1^3 + x_2^3 - x_1 + x_2 = 0,\ x_1 (x_1 - a_1) \geq 0 \right\} \cup \left\{(1, -1)\right\}, \\
{\rm WVAL} &=& \left\{(z_1, z_2) \in \mathbb{R}^2 \mid z_1 + z_2 = 0, z_1 (z_1 - 2) \geq 0 \right\} \cup \left\{(0, -2)\right\},
\end{eqnarray*}
are closed (see {\sc Figure~\ref{fig0609-2}}).

We now show that the estimate in the statement~\eqref{Item4} cannot, in general, be made uniform with respect to $\bar{x}$ on a compact set, compare \cite{Hu2024}. To this end, for each integer $k \ge 1,$ let $z^k := \left(-\frac{1}{k},\frac{1}{k}\right) \in \mathrm{WVAL}$ 
and choose $x_k \in f^{-1}(z_k) \cap \Omega.$ Certainly $z_k \to 0$ and the sequence $x_k$ is bounded. 
Let $\bar{y} := (1, -1) \in \Omega.$ We have $\Omega(x_k) \subset \Omega \setminus \{\bar{y}\},$ and so
\begin{eqnarray*}
\mathrm{dist}(\bar{y}, \Omega(x_k)) &\ge& \mathrm{dist}(\bar{y}, \Omega \setminus \{\bar{y}\}) \ > \ 0.
\end{eqnarray*}
Since $\psi(x_k) = 0$ and $\phi(x_k, \bar{y}) = \frac{1}{k} \to 0,$ it follows that for any compact set $X \subset \mathbb{R}^n,$ there are no constants $c > 0$ and $\alpha > 0$ such that
\begin{eqnarray*}
\phi(x, y) + \psi(x) +[\phi(x, y) + \psi(x)]^\alpha &\geq& c\, \mathrm{dist}(y, \Omega(x))
\end{eqnarray*}
for all $(x, y)\in X \times \Omega.$
}\end{example}

The following simple example shows that, in general, the term $\left[\max_{i = 1, \ldots, m} \big(f_i(y) - \bar{z}_i \big) \right]^\alpha$ in the statements~\eqref{Item5} and~\eqref{Item6} of Theorem~\ref{Thm52} cannot be removed.

\begin{example}{\rm
Let $\Omega := \mathbb{R}^2$ and consider the polynomial function
\begin{eqnarray*}
f \colon \mathbb{R}^2 \rightarrow \mathbb{R}, \quad (x_1, x_2) \mapsto x_1^2 + x_2^4.
\end{eqnarray*}
Then we can check that $f$ is convenient, nondegenerate at infinity and bounded from below.
Observe that $L(0) = \mathrm{ESOL} = \mathrm{WSOL} = \{(0, 0)\}$ and $\mathrm{EVAL} = \mathrm{WVAL} = \{0\}.$
Moreover, for each $z \in \mathbb{R},$
\begin{eqnarray*}
L(z)
&=& \begin{cases}
\emptyset & \textrm{ if } z < 0, \\
\left\{ (x_1, x_2) \in \mathbb{R}^2 \mid x_1^2 + x_2^4 \leq z\right\} & \textrm{ otherwise.}
\end{cases}
\end{eqnarray*}
Hence there is no constant $c > 0$ such that
\begin{eqnarray*}
f(y) &\ge& c\, \mathrm{dist}(y, \mathrm{WSOL}) \quad \textrm{ for all } \quad y \in \mathbb{R}^2,
\end{eqnarray*}
while for $c := \frac{1}{2}$ we have
\begin{eqnarray*}
f(y) + \left[f(y)\right]^{\frac{1}{4}} &\ge& c\, \mathrm{dist}(y, \mathrm{WSOL}) \ = \  c\, \mathrm{dist}(y, \mathrm{ESOL}) \quad \textrm{ for all } \quad y \in \mathbb{R}^2.
\end{eqnarray*}
for all $y \in \mathbb{R}^2.$
}\end{example}

The following example shows that, in general, the upper H\"older semicontinuity of the set-valued mapping
$$L \colon \mathbb{R}^{m} \rightrightarrows \mathbb{R}^n, \quad z \mapsto f^{-1}(z - \mathbb{R}^m_+) \cap \Omega,$$
cannot be strengthened to the H\"older continuity even when the assumptions of Theorem~\ref{Thm52} are satisfied.

\begin{example}\label{Example64}{\rm
Let $\Omega := \mathbb{R}^2$ and consider the polynomial function
\begin{eqnarray*}
f \colon \mathbb{R}^2 \to \mathbb{R}, && (x_1, x_2) \mapsto (x_1^2 + x_2^2-1)^2.
\end{eqnarray*}
Clearly, $f$ is bounded from below with ${\rm WVAL} = \{0\}$ and the optimal set $\mathrm{argmin}_{x \in \mathbb{R}^n} f(x)$ is just the unit circle in $\mathbb{R}^2.$
By definition, it is not hard to check that $f$ is convenient and nondegenerate at infinity. Moreover, we have (see {\sc Figure}~\ref{fig0518-1})
\begin{eqnarray*}
L(z)
\! &=& \! \begin{cases}
\emptyset & \text{if } z < 0, \\
\left\{(x_1,x_2) \mid 1 - \sqrt{z} \leq x_1^2 + x_2^2 \leq 1 + \sqrt{z}\right\} & \text{if } 0 \leq z < 1,\\
\left\{(x_1,x_2) \mid x_1^2 + x_2^2 \leq 1 + \sqrt{z}\right\} & \text{if } z \geq 1.
\end{cases}
\end{eqnarray*}
A direct computation shows that for any $\bar z \in \mathbb{R},$ there exist constants $c > 0$ and $r > 0$ such that
\begin{eqnarray*}
L(z) &\subset& L(\bar z) + c |z - \bar z|^{\alpha} \mathbb{B} \quad \textrm{ for all } \quad z \in \mathbb{B}_r(\bar z),
\end{eqnarray*}
where
\begin{eqnarray*}
\alpha &:=&
\begin{cases}
\frac{1}{2} & \textrm{ if } \bar{z} = 0, \\
1 & \textrm{ otherwise.}
 \end{cases}
\end{eqnarray*}
On the other hand, since $L(z') = \emptyset$ when $z' < 0,$ there are no constants $c > 0, \alpha > 0,$ and $r > 0$ such that
\begin{eqnarray*}
L(z) &\subset& L(z') + c |z - z'|^{\alpha} \mathbb{B} \quad \textrm{ for all } \quad z, z' \in \mathbb{B}_r(0).
\end{eqnarray*}
In particular, $L$ is not H\"older continuous.

\begin{figure}[h!]
	\centering
\begin{tikzpicture}[scale=2]
 \draw[thick,->] (0,0,0.7) -- (0,0,3) node[left,below] {$x_1$};
 \draw[thick,->] (1,0,0) -- (2.5,0,0) node[right] {$x_2$};
 \draw[dashed,->] (0,-0.2,0) -- (0,2.7,0) node[above] {$z$};
 \draw[dashed] (0,0,-3) -- (0,0,0.7) node[left,below] {$\quad$};
 \draw[dashed] (-1.5,0,0) -- (1,0,0) node[right] {$\quad$};

 \fill[blue] (2,0) circle (1.0pt) node[below] {$2$};
 \fill[blue] (1,0) circle (1.0pt) node[below] {$1$};
 \fill[blue] (-1,0) circle (1.0pt) node[below] {$-1$};
 \fill[red] (0,1) circle (1.0pt) node[above right] {$1$};
 \fill[gray] (0,0.55) circle (1.0pt) node[right] {$z$};
 \fill[blue] (0,0) circle (1.0pt) node[below right] {$0$};
 \fill[blue] (0,0,0.54) circle (1.0pt) node[below right] {$1$};
 \fill[blue] (0,0,1.08) circle (1.0pt) node[below right] {$2$};
 \fill[red] (0,1) circle (1.0pt) node[left] {$\quad$};
 \draw[purple,fill=white] (0,1.8) circle (1.0pt) node[right] {$z'$};

 \draw[dashed, gray] (0,0.56) ellipse (0.5 and 0.15);

 \draw[thick, blue] (-1,0) arc (180:360:1 and 0.2);
 \draw[dashed, blue] (1,0) arc (0:180:1 and 0.2);

 \draw[thick, gray] (-1.33,0.6) arc (180:360:1.33 and 0.28);
 \draw[dashed, gray] (1.33,0.6) arc (0:180:1.33 and 0.28);

 \draw[thick, red] (-1.41,1) arc (180:360:1.41 and 0.36);
 \draw[dashed, red] (1.41,1) arc (0:180:1.41 and 0.36);

 \draw[thick, purple] (-1.52,1.8) arc (180:360:1.52 and 0.41);
 \draw[dashed, purple] (1.52,1.8) arc (0:180:1.52 and 0.41);

 \node[thick, purple] at (1.9,1.81) {$f^{-1}(z')$};
 \node[thick, red] at (1.9,.9) {$f^{-1}(1)$};
 \node[thick, gray] at (1.9,0.5) {$f^{-1}(z)$};

 \draw[thick, black, domain=-1.6:-1, samples=100] plot (\x, {((\x)^2-1)^2});
 \draw[thick, black, domain=1:1.6, samples=100] plot (\x, {((\x)^2-1)^2});
 \draw[thick, black, dashed, domain=-1:1, samples=100] plot (\x, {((\x)^2-1)^2});
\end{tikzpicture}
\caption{\small The graph of $f(x_1,x_2) = (x_1^2 + x_2^2-1)^2$ and its level sets in Example~\ref{Example64}.}
\label{fig0518-1}
\end{figure}
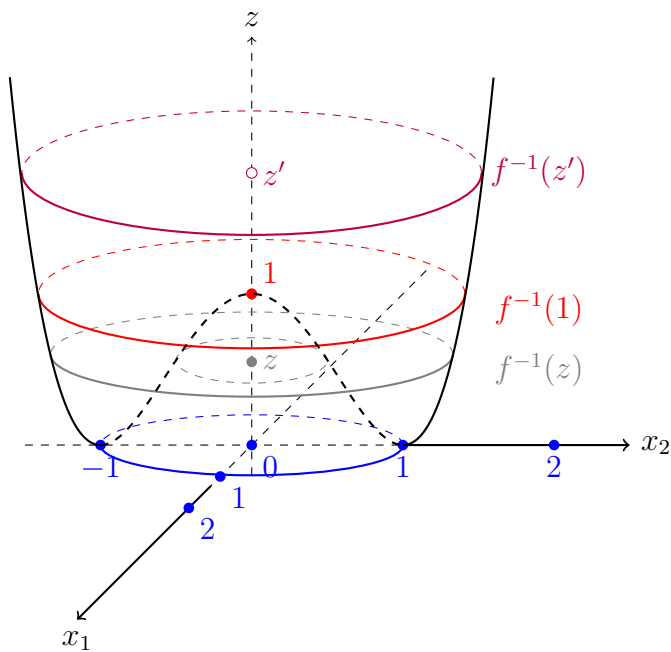
}\end{example}

In Theorem~\ref{Thm52}, the condition that $f$ is nondegenerate at infinity cannot be dropped as shown by the following example.

\begin{example}\label{Example65}{\rm
Let $\Omega := \mathbb{R}^2$ and consider the polynomial function
\begin{eqnarray*}
f \colon \mathbb{R}^2 \to \mathbb{R}, && (x_1, x_2) \mapsto (x_1 + x_2)^2 + ((x_1+x_2)x_2 - 1)^2.
\end{eqnarray*}
By definition, $\Gamma(f) = {\rm conv}\{(0,0),(2,0),(2,2),(0,4)\},$ and hence $f$ is convenient (see {\sc Figure~\ref{fig0519-0}}).
\begin{figure}[h!]
	\centering
\begin{tikzpicture}[scale=1.8]
  \draw[->] (-0.5,0) -- (2.5,0) node[right] {$\kappa_1$};
  \draw[->] (0,-0.5) -- (0,2.5) node[above] {$\kappa_2$};
  \fill[blue] (2,0) circle (1.0pt) node[below] {$4$};
  \fill[blue] (1.5,0) circle (1.0pt) node[below] {$3$};
  \fill[blue] (0,1.5) circle (1.0pt) node[left] {$3$};
  \fill[blue] (0,2) circle (1.0pt) node[left] {$4$};
  \fill[blue] (1,0) circle (1.0pt) node[below] {$2$};
  \fill[blue] (0.5,0) circle (1.0pt) node[below] {$1$};
  \fill[blue] (0,0.5) circle (1.0pt) node[left] {$1$};
  \fill[blue] (0,1) circle (1.0pt) node[left] {$2$};
  \fill[blue] (0,0) circle (1.0pt) node[below left] {$0$};
  \draw[thick, red] (1,1) -- (0,2);
  \draw[thick, red] (1,1) -- (1,0);
  \draw[thick, red] (0,0) -- (0,2);
  \draw[thick, red] (0,0) -- (1,0);

  \fill[red] (0,0) circle (1.2pt) node[right] {$\quad$};
  \fill[red] (0.5,1.5) circle (1.2pt) node[right] {$\quad$};
  \fill[red] (1,0) circle (1.2pt) node[right] {$\quad$};
  \fill[red] (0,2) circle (1.2pt) node[right] {$\quad$};
  \fill[red] (1,1) circle (1.2pt) node[right] {$\quad$};
  \fill[red] (0,1) circle (1.2pt) node[right] {$\quad$};
  \node at (0.5,0.8) {$\Gamma(f)$};
\end{tikzpicture}
\caption{\small The Newton polyhedron of $f$ in Example~\ref{Example65}.}
\label{fig0519-0}
\end{figure}
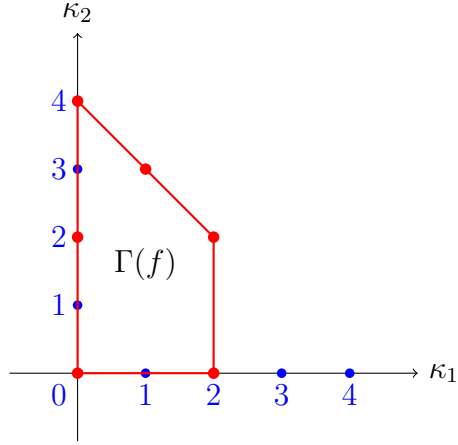
Now, let $I := \{1\},$ $J := \{1, 2\},$ and $q := (-1, -1) \in \mathbb{R}^2.$ We have $\min\{q_1, q_2\} = -1 < 0,$ and
\begin{itemize}
\item[$\circ$] $d = d(q, \Gamma(f)) = -4 < 0;$ and
\item[$\circ$] $\Delta = \Delta(q, \Gamma(f)) = {\rm conv}\{(2,2), (0,4)\}\subset \mathbb{R}^{J}.$
\end{itemize}
By definition, then
\begin{eqnarray*}
f_{\Delta}(x_1, x_2) &=& x_1^2x_2^2 + x_2^4 + 2 x_1x_2^3.
\end{eqnarray*}
Clearly, the point $(-1, 1) \in (\mathbb{R}^*)^J$ is a solution of the system
\begin{eqnarray*}
f_{\Delta} (x_1, x_2) \ = \ 0 \quad \textrm{ and } \quad  \nabla f_{\Delta}(x_1, x_2) & = & 0.
\end{eqnarray*}
Therefore, the polynomial $f$ is degenerate at infinity.
Observe that $f(\mathbb{R}^2) = (0, \infty),$ and hence the polynomial $f$ is bounded from below but does not attain its infimum.
}\end{example}

In Theorem~\ref{Thm52},  the condition that the mapping $f$ is convenient cannot be dropped as shown by the following examples.
\begin{example}\label{Example66}{\rm
Let $\Omega := \mathbb{R}^2$ and consider the polynomial mapping
\begin{eqnarray*}
f \colon \mathbb{R}^2 \to \mathbb{R}^2, && (x_1, x_2) \mapsto \left(x_1^2, (x_1x_2 - 1)^2 \right).
\end{eqnarray*}
By definition, $\Gamma(f_1) = \left\{(2,0)\right\}$ and $\Gamma(f_2) = {\rm conv}\left\{(0,0),(2,2)\right\},$ and hence $f$ is not convenient (see {\sc Figure}~\ref{fig0605-1}).
Moreover, it is easy to check that the mapping $f$ is nondegenerate at infinity.
On the other hand, it is not hard to see that
\begin{eqnarray*}
f(\mathbb{R}^2) &=& \left\{(z_1,z_2) \in \mathbb{R}^2 \mid z_1 > 0,\ z_2 \geq 0\right\} \cup \left\{(0, 1)\right\}
\end{eqnarray*}
(see {\sc Figure~\ref{fig0507-1}}(b)); in particular, the mapping $f$ has nonempty bounded sections and empty sections.
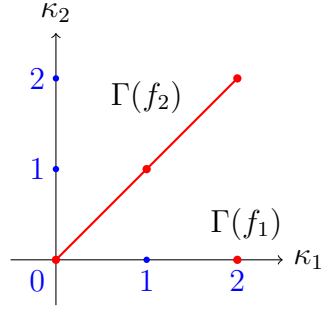
\begin{figure}[h!]
	\centering
\begin{tikzpicture}[scale=1.2]
  \draw[->] (-0.5,0) -- (2.5,0) node[right] {$\kappa_1$};
  \draw[->] (0,-0.5) -- (0,2.5) node[above] {$\kappa_2$};
  \fill[blue] (2,0) circle (1.0pt) node[below] {$2$};
  \fill[blue] (1,0) circle (1.0pt) node[below] {$1$};
  \fill[blue] (0,1) circle (1.0pt) node[left] {$1$};
  \fill[blue] (0,2) circle (1.0pt) node[left] {$2$};
  \fill[blue] (0,0) circle (1.0pt) node[below left] {$0$};
  \draw[thick, red] (0,0) -- (2,2);
  \fill[red] (0,0) circle (1.3pt) node[right] {$\quad$};
  \fill[red] (1,1) circle (1.3pt) node[right] {$\quad$};
  \fill[red] (2,2) circle (1.3pt) node[right] {$\quad$};
  \fill[red] (2,0) circle (1.3pt) node[right] {$\quad$};
  \node at (2.1,0.4) {$\Gamma(f_1)$};
  \node at (1,1.8) {$\Gamma(f_2)$};
\end{tikzpicture}
\caption{\small The Newton polyhedra of $f_1$ and $f_2$ in Example~\ref{Example66}.}
\label{fig0605-1}
\end{figure}

Now, a careful calculation shows that
\begin{eqnarray*}
{\rm ESOL} &=& \left\{(0, x_2)\mid x_2 \in \mathbb{R} \right\}, \\
{\rm WSOL} &=& \left\{(0, x_2)\mid x_2 \in \mathbb{R} \right\} \cup \left\{(x_1, x_2) \mid x_1 x_2 = 1 \right\}, \\
{\rm EVAL} &=& \{(0,1)\},
\end{eqnarray*}
which are all closed, while the set
\begin{eqnarray*}
{\rm WVAL} &=& \{(0,1)\}\cup\{(z_1, 0) \mid z_1 > 0\}
\end{eqnarray*}
is not closed (see {\sc Figure~\ref{fig0507-1}}).
\begin{figure}[h!]
\centering
\begin{tikzpicture}[scale=0.8]
\begin{scope}[shift={(-4,0)}]
  \draw[->] (-3.1,0) -- (3.1,0) node[right] {$x_1$};
  \draw[->] (0,-2.5) -- (0,4) node[above] {$x_2$};
  \draw[domain=0.4:3, samples=50, smooth, very thick, red]
    plot (\x, {1/\x});
  \draw[domain=-3:-0.4, samples=50, smooth, very thick, red]
    plot (\x, {1/\x});
 \draw[very thick, red] (0,-2.4) -- (0,3.2);
 \draw[yellow] (-0.02,-2.3) -- (-0.02,3.1);
  \node at (-0.8,1.5) {ESOL};
  \node at (1.4,-1.1) {WSOL};

  \coordinate (A) at (-0.8,1.1);
  \coordinate (B) at (-0.1,0.5);


  \draw[->, dashed] (A) to[bend right=60] (B);

\draw[->,dashed] (1.3,-0.8) -- (0.1,-0.1) node[right] {};
\draw[->,dashed] (1.4,-0.8) -- (1.6,0.5) node[right] {};
\draw[->,dashed] (0.75,-1.2) -- (-0.5,-1.6) node[right] {};
  \node at (0,-3.5) {${\rm (a)}$};
\end{scope}

\begin{scope}[shift={(4,0)}]
\draw[->] (-2.5,0) -- (4,0) node[right] {$z_1$};
\draw[->] (0,-2.5) -- (0,4) node[above] {$z_2$};

\fill[gray!30, opacity=0.6] (0.02,0) rectangle (3.7,3.7);

\draw[white, line width=2pt] (0,0) -- (0,3.8);
\draw[dashed, thin] (0,0) -- (0,3.8);
\draw[blue, thick] (0,0) circle (0.06);
\draw[blue, line width=1pt] (0.06,0) -- (3.7,0);

\fill[blue] (0,1) circle (2pt) node[above right] {$(0,1)$};

\node at (1.8,1.8) {$f(\mathbb{R}^2)$};

  \node at (0.5,-3.5) {${\rm (b)}$};

\end{scope}
\end{tikzpicture}
\caption{\small (a) The efficient solution set and weakly efficient solution set; \ (b) the image $f(\mathbb{R}^2)$ and the set of efficient values in Example~\ref{Example66}.}
\label{fig0507-1}
\end{figure}
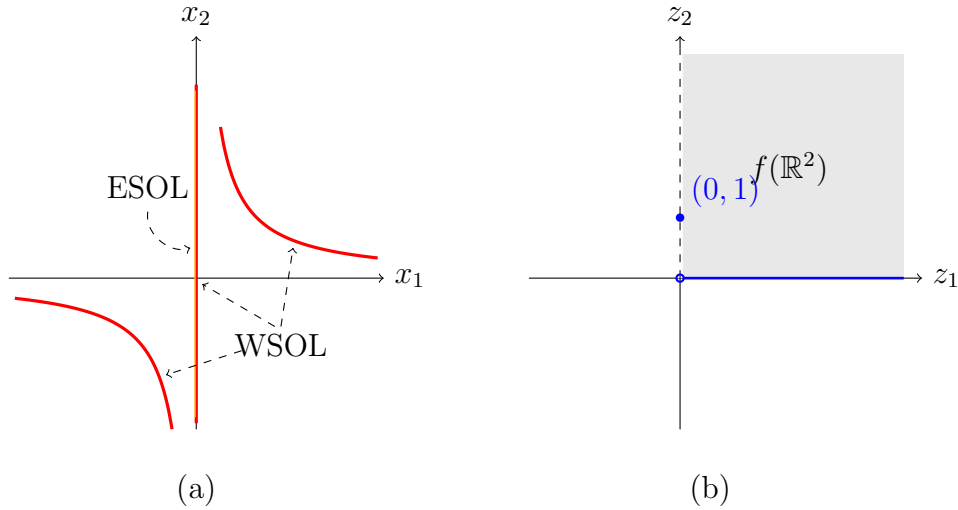
}\end{example}

\begin{example}\label{Example67}{\rm
Let $\Omega := \mathbb{R}^2$ and consider the polynomial mapping
\begin{eqnarray*}
f \colon \mathbb{R}^2 \to \mathbb{R}^2, && (x_1, x_2) \mapsto \left(x_1, x_2^2\right).
\end{eqnarray*}
By definition, $\Gamma(f_1) = \left\{(1,0)\right\}$ and $\Gamma(f_2) = \left\{(0,2)\right\},$ and hence $f$ is not convenient.
Moreover, it is easy see that the mapping $f$ is nondegenerate at infinity.
On the other hand, it is not hard to see that
\begin{eqnarray*}
f(\mathbb{R}^2) &=& \left\{(z_1, z_2) \in \mathbb{R}^2 \mid z_2 \geq 0 \right\}
\end{eqnarray*}
(see {\sc Figure~\ref{fig0512-3}}); in particular, the mapping $f$ has no nonempty bounded sections, but it has empty sections.
\begin{figure}[h!]
	\centering
\begin{tikzpicture}[scale=0.9]
  \draw[->] (-3.5,0) -- (3.5,0) node[right] {$z_1$};
  \draw[->] (0,-0.5) -- (0,3.5) node[above] {$z_2$};
  \fill[gray!30, opacity=0.5] (-3.4,0) rectangle (3.4,3.2);
  \draw[thick, blue] (-3.5,0) -- (3.5,0);
  \node at (1.5,1.5) {$f(\mathbb{R}^2)$};
\end{tikzpicture}
\caption{\small The image $f(\mathbb{R}^2)$ in Example~\ref{Example67}.}
\label{fig0512-3}
\end{figure}
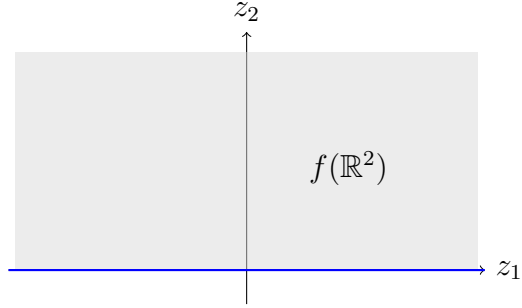

Now, a simple calculation shows that $\mathrm{ESOL} = \mathrm{EVAL} = \emptyset$ while
$\mathrm{WSOL} = \mathrm{WVAL} = \mathbb{R} \times \{0\} \ne \emptyset.$
}\end{example}

The following example shows that Corollary~\ref{Corollary51} can fail in some situations where the equation~\eqref{Eqn10} is replaced by the equivalent conditions (c1)--(c7) in Theorem~\ref{Thm52}.

\begin{example}\label{Example68}{\rm
Let $\Omega := \mathbb{R}^2_+$ and consider the polynomial mapping
\begin{eqnarray*}
f := (f_1, f_2) \colon \mathbb{R}^2 \to \mathbb{R}^2, && (x_1, x_2) \mapsto \left(-x_1^2 + x_2,\ x_1 + x_2\right).
\end{eqnarray*}
By definition, $\Gamma(f_1) = {\rm conv}\left\{(2,0),(0,1)\right\}$ and $\Gamma(f_2) = {\rm conv}\left\{(1,0),(0,1)\right\},$ and hence $f$ is convenient (see {\sc Figure~\ref{fig0821-1}}).
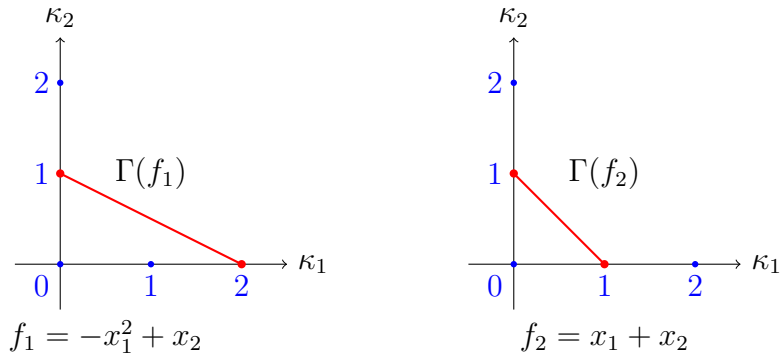
\begin{figure}[h!]
	\centering
\begin{tikzpicture}[scale=1.2]

\begin{scope}[shift={(-2.5,0)}]
  \draw[->] (-0.5,0) -- (2.5,0) node[right] {$\kappa_1$};
  \draw[->] (0,-0.5) -- (0,2.5) node[above] {$\kappa_2$};
  \fill[blue] (2,0) circle (1.0pt) node[below] {$2$};
  \fill[blue] (1,0) circle (1.0pt) node[below] {$1$};
  \fill[blue] (0,1) circle (1.0pt) node[left] {$1$};
  \fill[blue] (0,2) circle (1.0pt) node[left] {$2$};
  \fill[blue] (0,0) circle (1.0pt) node[below left] {$0$};
  \draw[thick, red] (2,0) -- (0,1);
  \fill[red] (2,0) circle (1.3pt) node[right] {$\quad$};
  \fill[red] (0,1) circle (1.3pt) node[right] {$\quad$};
  \node at (1,1) {$\Gamma(f_1)$};
  \node at (0.5,-0.8) {$f_1 = -x_1^2 + x_2$};
\end{scope}

\begin{scope}[shift={(2.5,0)}]
  \draw[->] (-0.5,0) -- (2.5,0) node[right] {$\kappa_1$};
  \draw[->] (0,-0.5) -- (0,2.5) node[above] {$\kappa_2$};
  \fill[blue] (2,0) circle (1.0pt) node[below] {$2$};
  \fill[blue] (0,2) circle (1.0pt) node[left] {$2$};
  \fill[blue] (1,0) circle (1.0pt) node[below] {$1$};
  \fill[blue] (0,1) circle (1.0pt) node[left] {$1$};
  \fill[blue] (0,0) circle (1.0pt) node[below left] {$0$};
  \draw[thick, red] (1,0) -- (0,1);
  \fill[red] (1,0) circle (1.3pt) node[right] {$\quad$};
  \fill[red] (0,1) circle (1.3pt) node[right] {$\quad$};
  \node at (1,1) {$\Gamma(f_2)$};
  \node at (1,-0.8) {$f_2 = x_1 + x_2$};
\end{scope}
\end{tikzpicture}
\caption{\small The Newton polyhedra of $f_1$ and $f_2$ in Example~\ref{Example68}.}
\label{fig0821-1}
\end{figure}
Moreover, by definition, we can check that the restriction $f|_{\Omega}$ is nondegenerate at infinity.

On the other hand, it is not hard to see that
\begin{eqnarray*}
f(\Omega) &=& \left\{(z_1, z_2) \in \mathbb{R}^2 \mid z_2 \geq 0,\ z_2 \geq z_1 \geq -z_2^2\right\}
\end{eqnarray*}
(see {\sc Figure~\ref{fig0821-2}}); in particular, the restriction $f|_\Omega$ has nonempty bounded sections. A direct calculation shows that
$(-1, 0) \in [f(\Omega)]^{\oplus}\cap(-\mathbb{R}^2_+).$ By Corollary~\ref{Corollary51}, the problem~\eqref{MP} does not have a properly efficient solution.

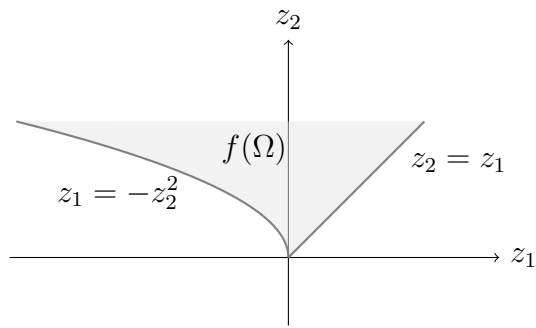
\begin{figure}[h!]
	\centering
\begin{tikzpicture}[scale=0.9]
  \draw[->] (-4.1,0) -- (3.1,0) node[right] {$z_1$};
  \draw[->] (0,-1) -- (0,3.2) node[above] {$z_2$};
  \fill[gray!20, opacity=0.5]
    (0,0) -- plot[domain=0:2, samples=100] ({-\x*\x},\x)
    -- (2,2) -- cycle;
  \draw[thick, gray] plot[domain=0:2, samples=100] (\x, \x)
    node[right] {};
  \draw[thick, gray] plot[domain=-2:0, samples=100] (-{\x*\x}, -\x)
    node[right] {};
  \node at (-0.5,1.6) {$f(\Omega)$};
  \node at (-2.5,0.95) {$z_1 = -z_2^2$};
  \node at (2.5,1.4) {$z_2 = z_1$};
\end{tikzpicture}
\caption{\small The image $f(\Omega)$ in Example~\ref{Example68}.}
\label{fig0821-2}
\end{figure}
}\end{example}

\subsection*{Data availability} There is no data included in this paper.

\subsection*{Declarations} No potential conflict of interest was reported by the authors.

\end{document}